\documentclass[preprint,3p,12pt]{elsarticle}
\usepackage{lineno,hyperref}
\usepackage{amsfonts}
\usepackage{amsmath,amssymb,epsfig}
\usepackage[utf8]{inputenc}
\usepackage[english]{babel}
\usepackage{mathtools}
\usepackage{helvet}         
\usepackage{courier}        
\usepackage{type1cm}
\usepackage{csquotes}
\usepackage{makeidx}        
\usepackage{graphicx}
\usepackage{multicol}
\usepackage{appendix}
\usepackage{stix}
\usepackage[bottom]{footmisc}
\usepackage[dvipsnames]{xcolor}
\usepackage{multicol}
\usepackage[bottom]{footmisc}     
\usepackage{hyperref} 
\usepackage{tikz}
\usepackage{float}  
\usepackage{helvet,overpic} 
\usepackage{subfig}
\usepackage{soul}
\usepackage{amsthm}

\newtheorem{definition}{Definition}

\usepackage{nomencl}
\usepackage[linesnumbered,ruled,vlined]{algorithm2e}

\SetCommentSty{mycommfont}
\SetKwInput{KwInput}{Input}                
\SetKwInput{KwOutput}{Output}
\makenomenclature
\hypersetup{pdfauthor={Kiumars Sharifmoghaddam}}
\usepackage[all]{xy}
\usepackage{amsthm}    
\usepackage{amsmath,amsfonts}                     
\usepackage{type1cm}             
\usepackage{makeidx}                      
\usepackage{multicol}          
\usepackage[bottom]{footmisc}
\usepackage{epsfig}     
\usepackage{psfrag,rotating}     
\usepackage{amssymb,floatflt,enumerate} 
\usepackage{amscd,psfrag} 
\usepackage[mathscr]{eucal}  
\usepackage{shadethm}           
\usepackage{overpic,contour}       
\contourlength{0.3mm}

\DeclareMathOperator{\sign}{sign}

\newtheorem{theorem}{\bf Theorem}[section]

\newtheorem{remark}{\bf Remark}[section]
\newtheorem{lemma}{\bf Lemma}[section]

\def\RR{{\mathbb R}}

\def\dach#1{\widehat{#1}}

\def\Vkt#1{{\mathbf #1}}

\newcommand{\mVkt}[1]{\dach{\Vkt #1}}

\begin{document}

\begin{frontmatter}

\title{Generalizing rigid-foldable tubular structures of T-hedral type}

\author{K. Sharifmoghaddam$^1$, R. Maleczek$^2$, G. Nawratil$^1$}
\address{$^1$Institute of Discrete Mathematics and Geometry \&  
	Center for Geometry and Computational Design, TU Wien \\
	$^2$i.sd - structure and design, Department of Design,
 Faculty of Architecture, Universität Innsbruck}

\address[mymainaddress]{Wiedner Hauptstrasse 8-10/104, Vienna 1040, Austria}

\cortext[mycorrespondingauthor]{Corresponding author: Kiumars Sharifmoghaddam ({ksharif@geometrie.tuwien.ac.at})}

\begin{abstract}
We introduce an alternative way of constructing  
continuous flexible
tubes and tubular structures based on a discrete, semi-discrete and smooth construction of surfaces known as 
 T-hedra  in the discrete case and 
profile-affine surfaces in the smooth setting, respectively. 
The geometric understanding of this method enables us to generalize discrete tubes with a rigid-foldability  and to extend the construction to smooth and semi-discrete tubes with an isometric deformation.
This achievement implies a unified treatment of continuous flexible structures, like surfaces and metamaterials,  composed of tubes and it is the base for a deeper study of zipper tubes, 
and their generalization. 
Moreover, we discuss a potential application of the presented structures for the design of foldable bridges.
\end{abstract}

\end{frontmatter}

\section{Introduction}\label{sec:Intro}

 Rigid-foldable tubular structures are lightweight and robust structures with special kinematic properties, which have been of interest to designers, architects and engineers over the past decade. Having one-degree-of-freedom (1-DoF) movement and other possible geometrical properties, such as flat-foldability, make them suitable to be utilized in kinetic architecture, transportable structures and, in general, transformable design.

\begin{figure}
\centering
\begin{overpic}[scale=0.25]{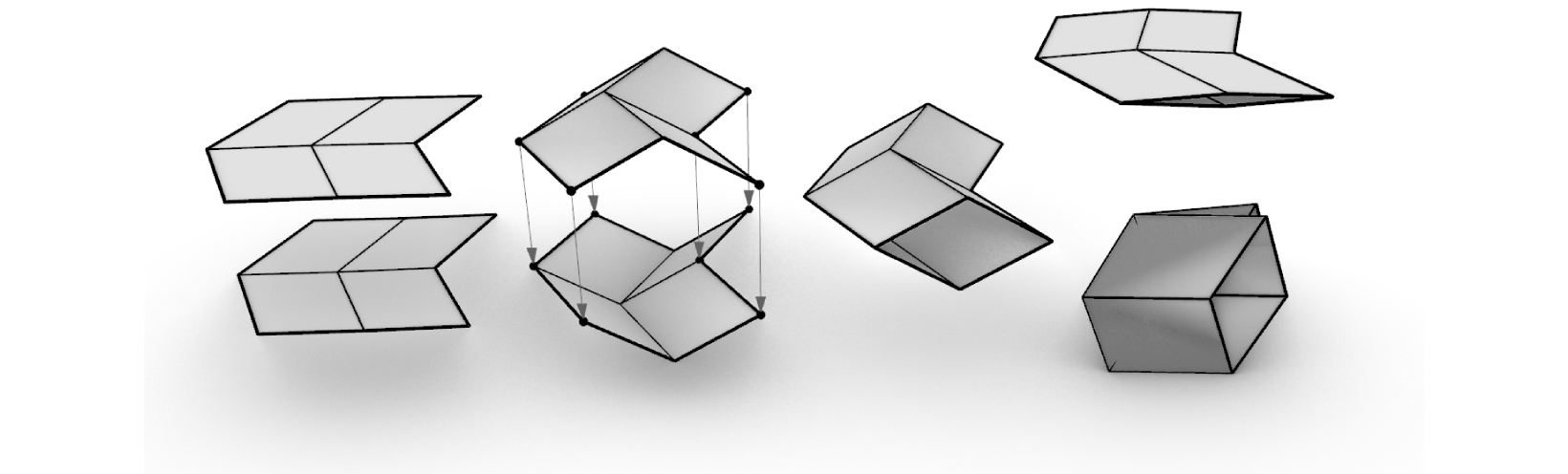}
{\small
\put(22,3){(a)} 
\put(41,3){(b)}
\put(59,3){(c)}
\put(75,3){(d)}
\put(75,20){(e)}
}
\end{overpic}
    \caption{Construction of a regular unit of a rigid-foldable tube with a rhombic cross-section. (a) Two congruent pieces of a regular Miura-ori unit in the developed state  and (b) in a half-folded one, respectively. (c) The tubular unit  which is (e) about to be flat-folded into its base plane and (d) into the carrier plane of the cross-section.}
    \label{One-Cell}
\end{figure}
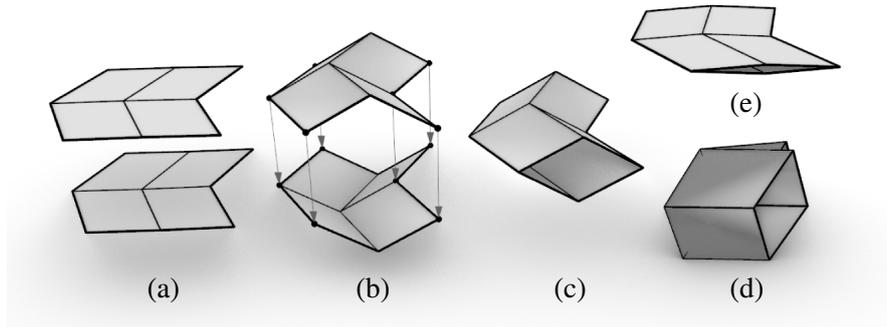

The simplest way to construct a unit of a rigid-foldable tube is to use two identical unit pieces of regular Miura-ori (consisting of four congruent parallelograms in a two-by-two pattern) and connect them in a mirrored way (Fig.\ref{One-Cell}). The result is a cylindrical structure that acts as a 1-DoF mechanism, which can be folded bidirectionally into two flat states. Furthermore, it is stiff in the third orthogonal direction, which is an important mechanical property and of interest for many potential applications. To be able to develop a wider design space, we take another geometrical approach to construct basic tubular units with a rigid-foldability. In our method, we take a cross-section (profile) and extrude it (or apply more complex transforms) along a path (trajectory) which is a more natural way of constructing a tube as such. Then, we generalize these initial inputs in many aspects, as well as the structures composed of the tubes, in the paper at hand. 

\subsection{Review (tubular structures)}\label{subsec:Review}

The geometric study of 1-DoF rigid-foldable tubular (or cylindrical as it was called back then) structures with rigid quadrilateral faces goes back to a paper by Tachi \cite{Tachi2009} in which the author proposed flat-foldable structures consisting of parallelograms where the cross-section of the tubes is generalized from a rhombus to polygons with two-folded rotational symmetries. The tube was the result of extruding the polygonal section by a regular zigzag trajectory as the simplest way to achieve flat-foldability.

Miura and Tachi in \cite{miura_synthesis} presented a new method of making tubes by connecting two identical layers of modified and trimmed Miura-ori pattern. These tubes, which can also be arranged to space-filling cellular structures, contain self-intersecting quads, resulting from the trimming process. Moreover, there is a fundamental difference to the construction presented in the paper at hand, which will be explained at the very end (Section \ref{sec:future}).

In \cite{tachi_freeform}, Tachi introduces a computational technique to obtain cylindrical rigid-foldable BDFFPQ (bidirectionally flat-foldable planar quadrilateral) meshes as a hybrid of generalized Miura-ori and eggbox patterns. The extrusion trajectory of the cross-sections were generalized zigzags in space. Later on, in \cite{Tachi2011} and \cite{tachi_transformables_2013}, Tachi introduced two new design methods to obtain flat-foldable tubes constructed along given spatial curves, which were approximated by a zigzag or a polyline in their initial state, respectively. Although structures composed of tubular cells already appeared in \cite{miura_synthesis}, the variety of possible connections of the constructed tubes had been investigated in more detail by the same authors in  \cite{tachi_iass2012}. Additionally, in \cite{tachi_transformables_2013}, a multilayered flat-foldable vault was shown as an example of structures composed of such tubes.

Another approach used by Schenk and Guest \cite{guest2013} was to obtain metamaterials by stacking Miura-ori layers. They connected two non-congruent Miura-ori patterns to make repetitive deltoid (kite shape) cross-sections and analyzed the stiffness of the structure. With a similar approach, Klett and Middendorf \cite{klett_iass2016} constructed sandwich structures and metamaterials using other variants of Miura-ori (e.g.\ curve folded) tessellated layers instead of usual edge-connected tubes.  

Theoretical development of geometric and kinematic methods, providing broader design space \cite{Chen,liu2016,davood} and  the identification of potential applications drew the attention of engineers to analyze the structural behaviour and to investigate the materialization of rigid-foldable tubular structures \cite{yasuda2013,cheung2014,Gattas2015,ma2018,you2021}.

Amongst other researchers, Filipov, Tachi and Paulino started to explore these structures from all the aspects mentioned above and published a series of papers, namely, on stiffness and flexibility optimization of a single tube \cite{FTPorigami6}, reconfigurable tubular structures and metamaterials \cite{pnas}, geometric investigations on the cross-sections \cite{FPT2015} and exploring the variety of couplings to obtain sandwich surfaces guided by a generating doubly curved surface \cite{tachietal2015} and with out-of-plane stiffness \cite{Filipovetal2019}. For the structural characteristics of such surfaces, see the work \cite{ishizawa} by Ishizawa and Tachi.    

\subsection{Overview}\label{subsec:Overview}

In this paper, we introduce an alternative way of constructing rigid-foldable tubes and tubular structures based on a discrete, semi-discrete and smooth construction of surfaces known as T-hedra (Section \ref{sec:T-hedra}) in the discrete case and profile-affine surfaces (Section \ref{sec:semitube}) in the smooth setting, respectively. The geometric understanding of this method enables us to generalize discrete tubes with a rigid-foldability (Section \ref{sec:general}) and to extend the construction to smooth and semi-discrete tubes with an isometric deformation (Section \ref{sec:semitube}). This achievement implies a unified treatment of continuous flexible structures composed of tubes (Section \ref{sec:Structure}), and it is the base for a deeper study of zipper tubes and their generalization (Section \ref{sec:Zipper}). 
In Section \ref{sec:bridge}, we discuss a potential application of the presented structures as a foldable bridge, give an outlook to future research and close with a list of open problems.
\section{Trapezoidal quad-surfaces (T-hedra)}\label{sec:T-hedra}

In mathematical terminology, surface is a general term. It can refer to smooth or discrete (also called polyhedral) surfaces, where the latter consists of vertices, edges and polygonal faces. Triangulation is the easiest way to discretize a surface. Since every three points in space are co-planar, a triangulated surface can always be constructed out of flat panels. If they are jointed by rotational hinges, one ends up with a transformable structure possessing several DoFs which are hard to control. 
In contrast, generic discrete surfaces composed of planar quads (PQ) connected by rotational joints in the combinatorics of the square grid are rigid. However, there exist special geometries allowing a 1-DoF flexibility. 

In this paper, we interchangeably use the wordings {\it continuous flexible},
{\it isometric deformable} and {\it rigid-foldable} for discrete surfaces, which change their spatial shape by rotating the hinges only and without any elastic or plastic deformations of panels. Therefore the 1-DoF motion of such PQ-surfaces can be controlled  by a single dihedral angle between hinged panels.   

The necessary and sufficient condition for a PQ-surface to be rigid-foldable is that  every $(3\times 3)$ complex has this property \cite{Schief}, which is an intrinsic one as it only relies on the corner angles of the quads. A partial classification of these building blocks was given by Stachel \cite{Stachel} and Nawratil \cite{naw1,naw2,naw3}, and a complete one was presented by Izmestiev \cite{Ivan}, containing more than 20 cases.  

The class, which contains the most types of rigid-foldable tubes reviewed in \ref{subsec:Review}, is called T-hedra. This name is a translation of the German term {\it T-Flache}, which was introduced by Graf and Sauer  \cite{graf} in 1931 for the first time. An important property of T-hedra, in contrast to general flexible quad-surfaces, is that they allow direct access to their spatial shape, despite the above-mentioned intrinsity of the rigid-foldability. The explicit geometric construction of T-hedra and the reasoning for their continuous flexion, were already known facts in German-language literature for several decades due to the work \cite{graf,sauer} of Graf and Sauer. Their basic ideas were summarized in the English written article \cite{SNRT2021}, which also presents explicit formulas for their isometric deformation (see also \cite{RT2021}). 
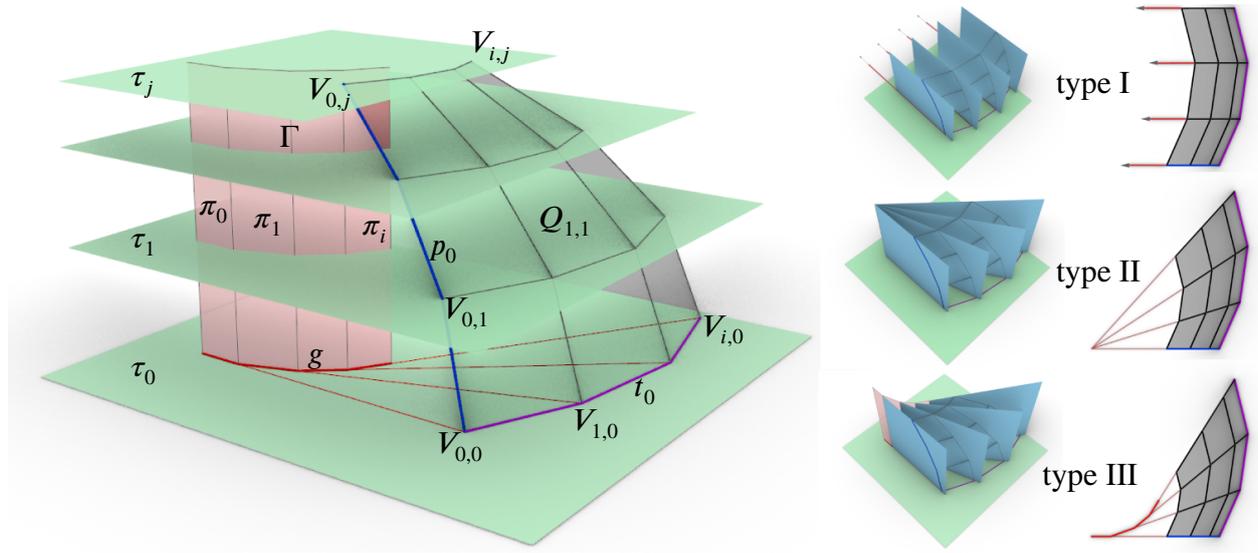
\begin{figure}
\centering
\begin{overpic}[scale=1.29]{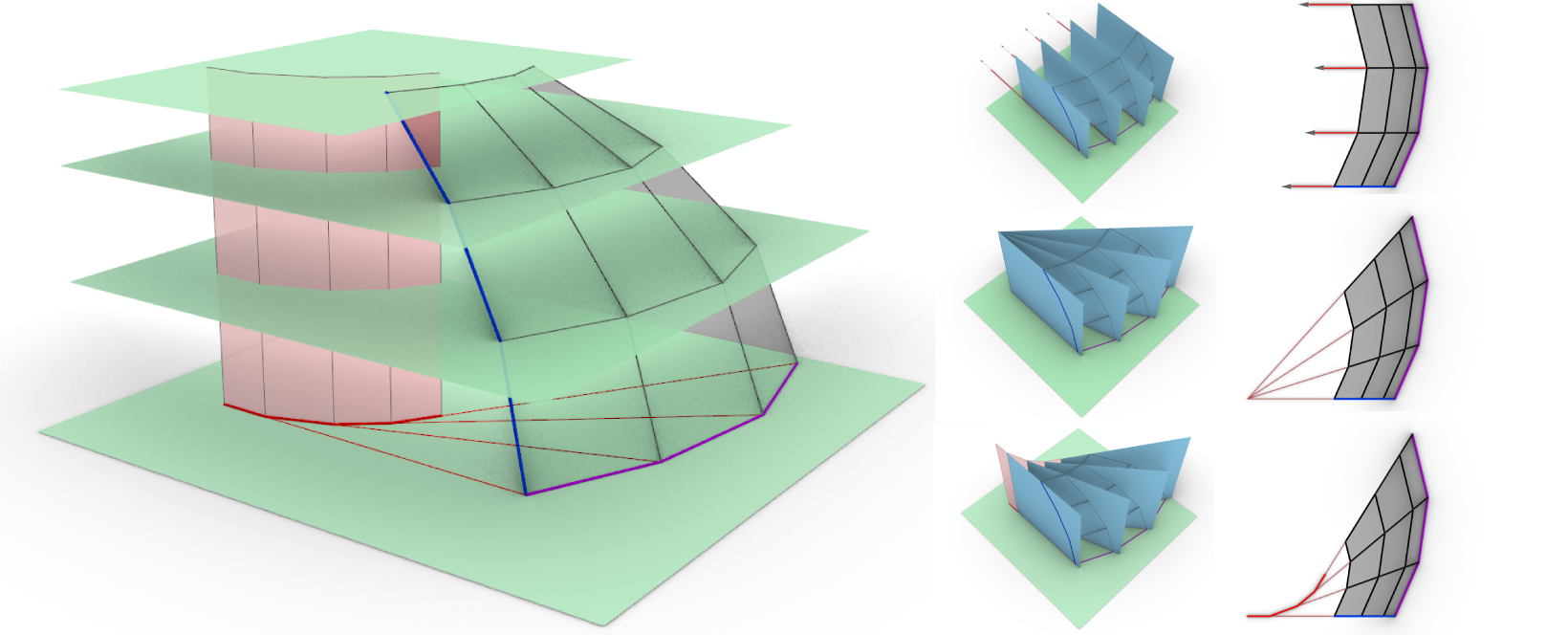}
{
\put(77,34){type I} 
\put(77,20){type II}
\put(76,5){type III}
\put(14,25){$\pi_0$}
\put(18,24){$\pi_1$}
\put(26,23.5){$\pi_i$}
\put(9,13){$\tau_0$}
\put(9,22.5){$\tau_1$}
\put(9,34.5){$\tau_j$}
\put(31.5,7.3){$V_{0,0}$}
\put(41.5,9){$V_{1,0}$}
\put(51,16){$V_{i,0}$}
\put(32,17.3){$V_{0,1}$}
\put(22,33.5){$V_{0,j}$}
\put(34,37){$V_{i,j}$}
\put(31,22){$p_0$}
\put(46,11.5){$t_0$}
\put(22,14){$g$}
\put(39,24){$Q_{1,1}$}
\put(20,30){$\Gamma$}
}
\end{overpic}
    \caption{(left) Geometric construction of a T-hedron. Illustrated with colors: Trajectory planes $\tau_0,\tau_1, \dots, \tau_j$ in green, profile polyline $p_0$ in dark blue, trajectory polyline $t_0$ in purple, guiding prism $\Gamma$ in light red and prism polyline $g$ in red. (right) Classification of T-hedra based on their top views, where the profile planes are shown in blue.}
    \label{P-T-planes}

\end{figure}
\subsection{Definitions}\label{subsec:Fund}

To construct a T-hedron, we start with two planar boundary polygons $p_0$ and $t_0$ located on two orthogonal planes $\pi_0$ and $\tau_0$ respectively, where the latter one is also known as {\it base plane}.
Moreover, $p_0$ and $t_0$ have their first vertex $V_{0,0}$ in common (see Fig.\ref{P-T-planes}, left). The remaining vertices of $t_0$ and $p_0$ are denoted by $V_{i,0}$ and $V_{0,j}$, respectively, with $i=1,...,I-1$ and $j=1,...,J-1$. Through each vertex $V_{i,0}$ of $t_0$ a plane $\pi_i$ passes, which is orthogonal to $\tau_0$ and will carry $p_i$. The envelope of these $\pi_i$ planes is called guiding prism $\Gamma$ and its intersection  with the plane $\tau_0$ yields the prism polygon $g$.

Starting from $p_0$, in each step of the construction, we obtain the planar polygon $p_i\in\pi_i$ via a parallel projection of $p_{i-1}$ into the plane $\pi_i$. The direction of this projection is determined by polygon edge $V_{i-1,0}V_{i,0}$ of $t_0$, which maps the vertices $V_{i-1,j}$ to $V_{i,j}$. By iterating this process (for $i=1,..,I-1$) we generate the vertices of the quad-mesh. Therefore the carrier planes $\tau_j$ of the polygons $t_j$ with the vertices $V_{0,j},V_{1,j},...,V_{I-1,j}$ are parallel to $\tau_0$. 
Moreover, each quad $Q_{i,j}$ is a planar trapezoidal face, as the sides $V_{i-1,j-1} V_{i,j-1}$ and $V_{i-1,j} V_{i,j}$ are parallel. The letter T, stands for trapezoidal, in the nomenclature T-hedron.

Without loss of generality we can assume $\pi_0 = XZ$-plane and $\tau_0 = XY$-plane, which results the corner vertex $V_{0,0}$ to lie on the $X$-axis.

\subsection{Subclassification of T-hedra based on their kinematic generation }\label{subsec:subclass}

In each step of the construction, the polygon $p_i$ can be generated kinematically from $p_{i-1}$ by a 
\begin{itemize}
     \item Type I: Translation
    \item Type II: Composition of a rotation and an axial dilatation with a stationary axis $l$
    \item Type III: Composition of a rotation and an axial dilatation with axis $\pi_{i-1} \cap \pi_i$ of $\Gamma$
\end{itemize}

These different types can easily be distinguished in the top view (Fig. \ref{P-T-planes}).
If we apply only pure rotations in cases II and III, we end up with discrete versions of rotation surfaces and molding surfaces, respectively. In both cases, all trapezoids are isosceles. For the discrete translation surfaces (type I), all trapezoids degenerate into parallelograms.
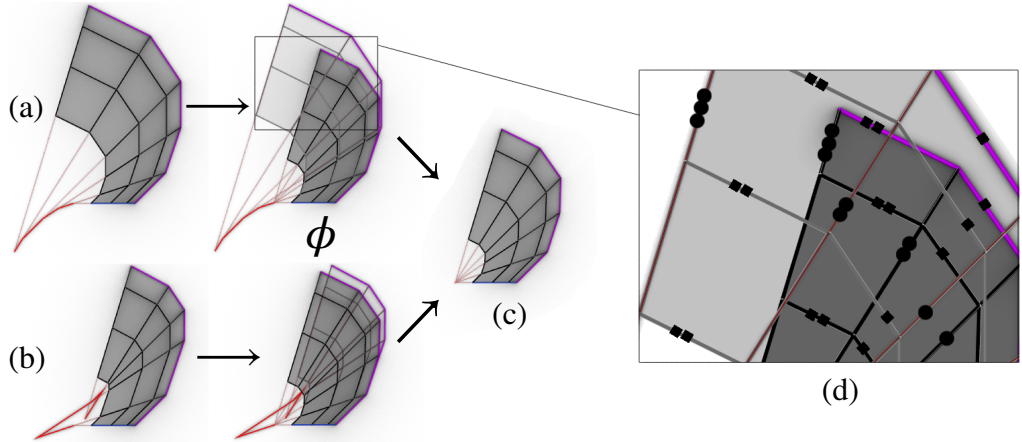
\begin{figure}
\centering
\begin{overpic}[scale=0.25]{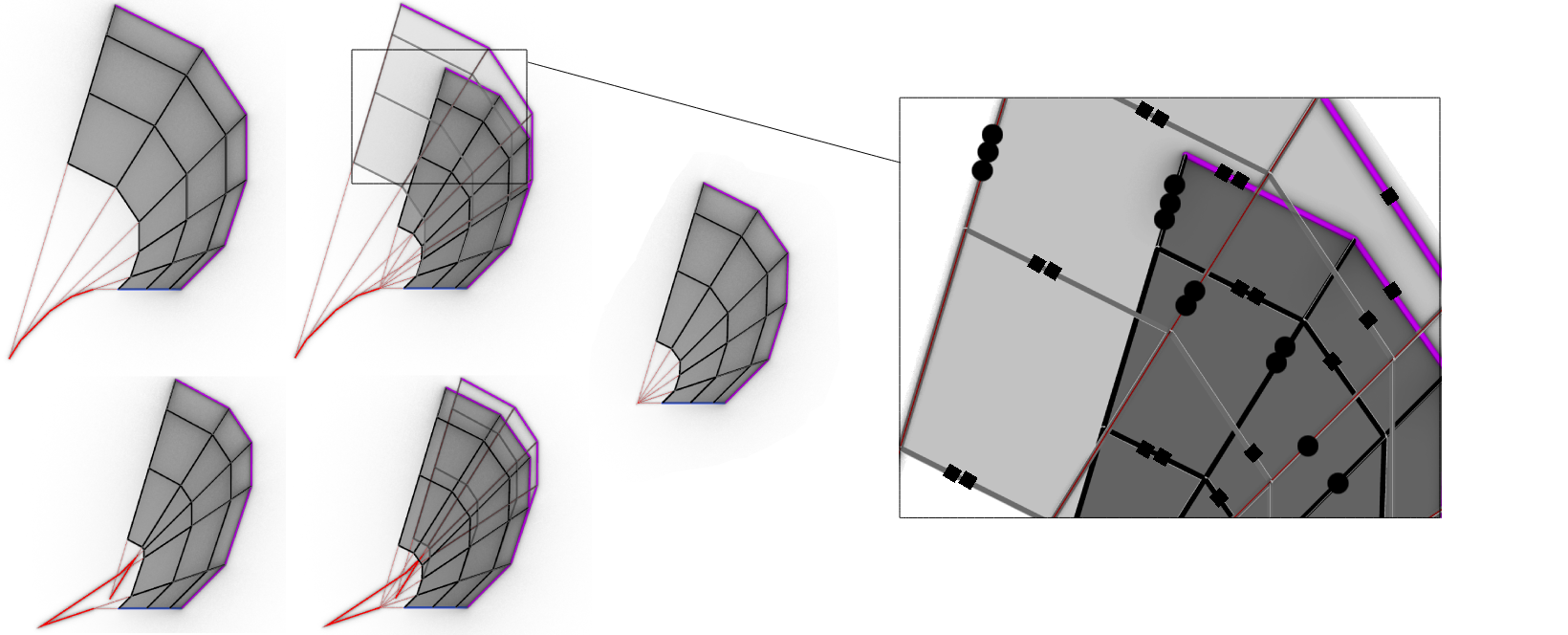}
{
\put(0,30){(a)}
\put(0,7){(b)}
\put(44,11){(c)}
\put(74,4){(d)}
\fontsize{16}{16}
\put(16,30){$\longrightarrow$}
\put(17,7){$\longrightarrow$}
\fontsize{18}{18}
\put(27,18){$\phi$}
\fontsize{20}{20}
\put(35,10){$\nearrow$}
\put(35,25){$\searrow$}

}
\end{overpic}
    \caption{The mapping $\phi$, which maps type III T-hedra to type II T-hedra. (a,b) Two different T-hedra of type III. (c) A T-hedron of type II is obtained by applying the mapping $\phi$ to (a) and (b), respectively. (d) A magnified frame of the mapping $\phi$ in which corresponding parallel profile edges are marked with round black dots and trajectory edges with black squares.}
    \label{mapping}

\end{figure}


\subsubsection{Relation between type II and type III}\label{sec:relation}
There is a mapping $\phi$ which maps type III T-hedra to type II T-hedra in the following way: All prism planes are translated such that they 
belong to a pencil of planes and the corresponding edges of the two T-hedra 
are parallel (see Fig.\ref{mapping}). The resulting type II T-hedron is uniquely defined up to translations, but the map $\phi$ is not injective. 
Therefore the deformation of a T-hedron of type III can always be reduced to the deformation of a T-hedron of type II (cf.\ \cite{graf, sauer, SNRT2021,RT2021}). 


\begin{figure}
\centering
\begin{overpic}[scale=0.25]{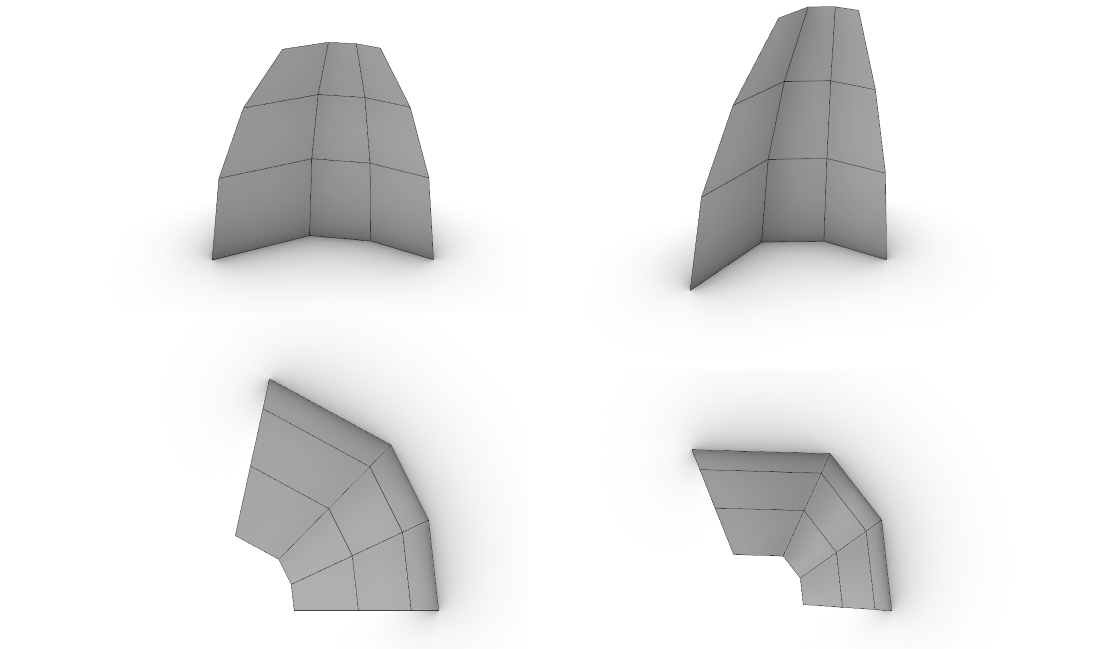}
{\small
\put(41,36){deformation}
\put(41,6){deformation}
\put(71,25){top view}
\put(31,25){top view}
\fontsize{14}{14}
\put(41,41){$T$}
\put(58,41){$T^*$}
\put(41,11){$T'$}
\put(57,11){$T'^* $}
\fontsize{22}{22}
\put(67,24){$\downarrow$}
\put(27,24){$\downarrow$}
\put(45,10){$\longrightarrow$}
\put(45,40){$\longrightarrow$}
}

\end{overpic}
    \caption{ A T-hedron $T$ and an  isometrically deformed version $T^*$ of $T$. The top views of $T$ and $T^*$ are denoted by $T'$ and $T'^*$,  respectively.}
    \label{3by3mesh-all4}

\end{figure}
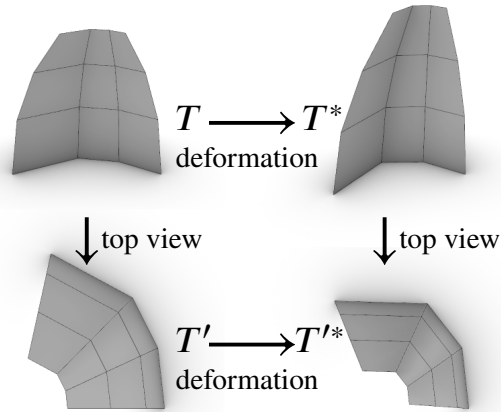
\subsection{Isometric deformation}\label{subsec:Iso}

As preparatory work for the formulation of the latter given Theorem \ref{closed-discrete-profile}, we sketch roughly the geometric proof for the isometric deformation of T-hedra according to Sauer and Graf \cite{graf}:

Due to the special geometry of a T-hedron $T$, the rigid-foldability of the structure can be reasoned by 
considering the top view $T'$
 (i.e. the  orthogonal  projection of $T$ into $XY$ plane). 
There exists a continuous 1-parametric deformation of $T'$ such that the deformed figure $T'^*$ can be 
seen as the top view of another T-hedron $T^*$, which originates from $T$ by a continuous isometric deformation (see Fig. \ref{3by3mesh-all4}).

One considers the deformation of the first strip of $T'$ by an affine transformation orthogonal to the 
parallel sides of the trapezoids, which yields the first strip of $T'^*$. The related distortion factor $t$ of this transformation 
can be chosen arbitrarily within some limits. The distortion factor of the next strip is already uniquely defined by the condition that 
the first and the second strip have to match. This procedure can be iterated, which finally yields $T'^*$.
Then it can easily be seen (cf.\ \cite{graf,SNRT2021}) that one can construct a surface $T^*$ over $T'^*$, which has the same inner geometry as $T$.

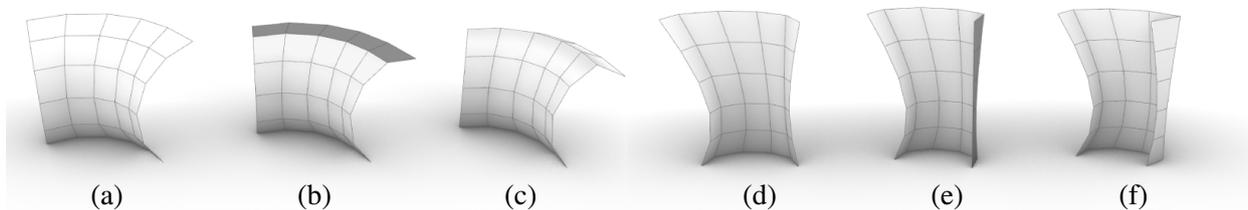
\begin{figure}[b]
\centering
\begin{overpic}[width=165mm]{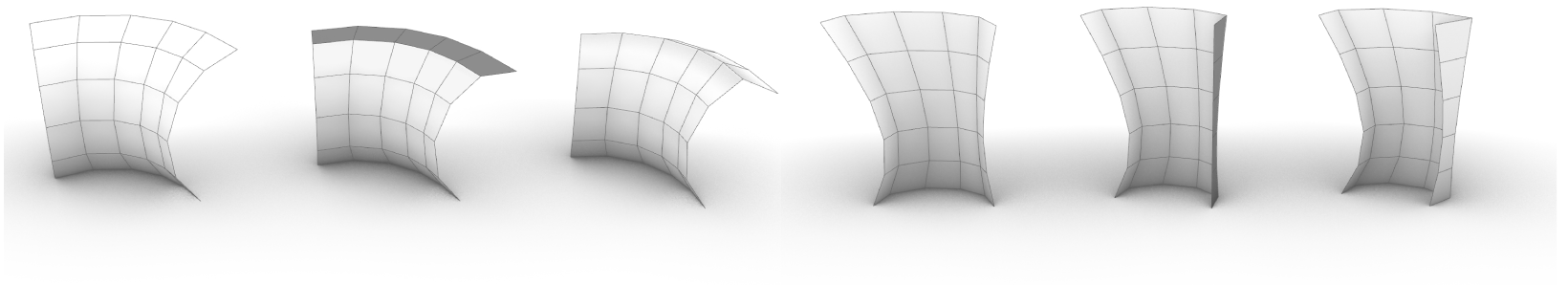}
{\small

\put(7,2){(a)} 
\put(23.5,2){(b)}
\put(40,2){(c)} 
\put(59,2){(d)} 
\put(74,2){(e)} 
\put(89,2){(f)} 
}
\end{overpic}
    \caption{ A T-hedron in (b) horizontal and (e) vertical bifurcation configurations (= horizontal and vertical
flexion limits). The two possible states for the top PQ-row and  the right PQ-column, which belong to different modes, are illustrated in (a,c) and (d,f), respectively.}
    \label{bifurcation}

\end{figure}

According to \cite{graf}, the flexion limits of a T-hedron are reached as soon as either a complete PQ-row gets parallel to the base plane (Fig.\ref{bifurcation}(b)), or a PQ-column (Fig.\ref{bifurcation}(e)) gets completely flat.  
\begin{remark}\label{discrete_bifurcation}
In the flexion limits the structure is also in a so-called bifurcation configuration, as it can transition from one working mode to another. As a consequence the T-hedron can only flex back if such a configuration is reached, but
 the PQ-strip can flip to the opposite side (Fig.\ref{bifurcation}(c,f)). This property for tubes is discussed latter on (cf.\ Section \ref{subsec:switch}) under the notion of {\it switching}. \hfill $\diamond$
\end{remark}

\section{Generalization of discrete tubes}\label{sec:general}

A close look at the literature on rigid-flexible tubes, composed of planar quads where each interior vertex has valence four, reveals that most of the known examples are T-hedra with a closed profile polyline. 
The only exceptions known to the authors are as follows:
\begin{enumerate}
    \item 
     Tubular structures with parallelogram cross-sections, which allow a more extensive design space of rigid-foldable tubes\footnote{\label{octahedron}Two neighboring tubular segments can be seen as a flexible octahedron with one pair of opposite points located on the plane at infinity \cite{N2015}. This circumstance also reasons the rigid-foldability of the referred tube structures.} than the one obtained by the T-hedral construction, as the trajectory polyline has not to be planar anymore (cf.\ \cite[Fig.\ 7]{FPT2015} and \cite{Tachi2011}). It can even be an arbitrary smooth space curve yielding a semi-discrete tube
    \cite{tachi_transformables_2013}. 
    Moreover, the cross-section can even be generalized according to the idea illustrated in \cite[Fig.\ 06]{tachi_transformables_2013}; i.e.
    the cross-section has not to be a parallelogram but a closed polyline where all edges are parallel to one of two given directions. Also see \cite[Fig.\ 8]{tachi_freeform}.
    \item
    Horn structures composed of Bricard octahedra of the third type \cite{tachi_horn} are examples of V-hedral\footnote{Graf and Sauer  \cite{graf} introduced discrete Voss surfaces under the name of V-hedra ({\it V-Flache} in German).} tubes.    
\end{enumerate}

In the following section, we unify and generalize the construction of tubes given in the literature (cf.\ Section \ref{subsec:Review}), which we have identified as being T-hedral.
\begin{figure}
\centering
\begin{overpic}[scale=0.8]{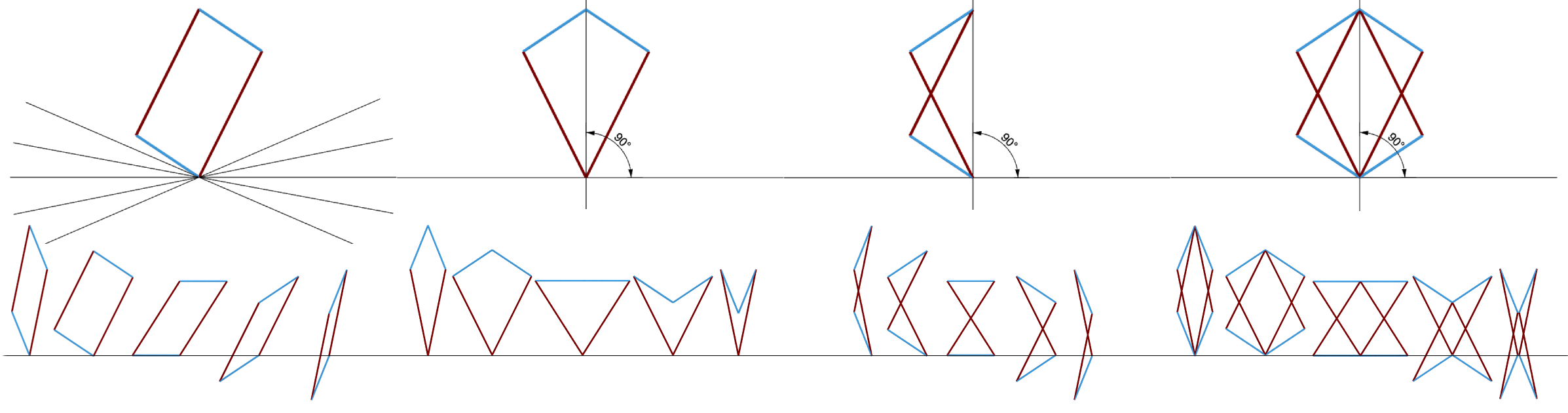}
{\small

}
\end{overpic}
    \caption{Types of possible quadrilateral cross-sections; from left to right: parallelogram, deltoid and anti-parallelogram. On the outmost right figure a cross-section is given which can be seen  as the union of 
    either two parallelgrams or deltoids or anti-parallelograms.}
    \label{all-quad-cross-sections}

\end{figure}
\subsection{Polygonal cross-section}\label{subsec:Poly}

Let us consider a closed planar polyline in the $XZ$-plane with starting point $S$ and endpoint $E=S$. This polyline  is suitable as a cross-section of a T-hedral tube if it 
remains closed under the deformation; i.e. $E^*=S^*$ has to hold. Clearly, if $E'^*=S'^*$ and  $E''^*=S''^*$ hold, where $''$ indicates the front view  (i.e. the  orthogonal projection into $YZ$plane), 
then this implies $E^*=S^*$. Due to the above-described principle of deformation, we get the property $E'^*=S'^*$ for free. Therefore we only have to ensure that the condition $E''^*=S''^*$ holds. This can be verified by the following procedure:

\begin{theorem}\label{closed-discrete-profile}
A closed cross-sectional polyline $P$ 
is compatible with the isometric deformations of T-type if  the following conditions hold:
\begin{enumerate}
\item 
Line-segments which have the same absolute value of the slope, belong to one class. The slopes can be given within the interval $[-90^{\circ};90^{\circ}]$, as 
every line-segment can be translated in a way that one of its endpoints is located in the origin and the other end-point is located in the first or fourth quadrant.
\item
For each class, we add up the oriented lengths of the individual line-segments, where the orientation can be defined as follows: 
We start at $S$ and run counter-clockwise through the polyline until we reach the point $E=S$. Line-segments pointing upwards 
are positive, and those pointing downwards are negative. We can assume without loss of generality that there exists a state of 
deformation where no line-segment is horizontal as otherwise, the structure 
cannot flex\footnote{Only within a flexion limit line-segments can be horizontal.}. 
\item
For each class this sum has to be zero\footnote{As the height changes during the deformation have to cancel out.}. 
\end{enumerate}
\end{theorem}
\begin{figure}[t]
\centering
\begin{overpic}[scale=0.68]{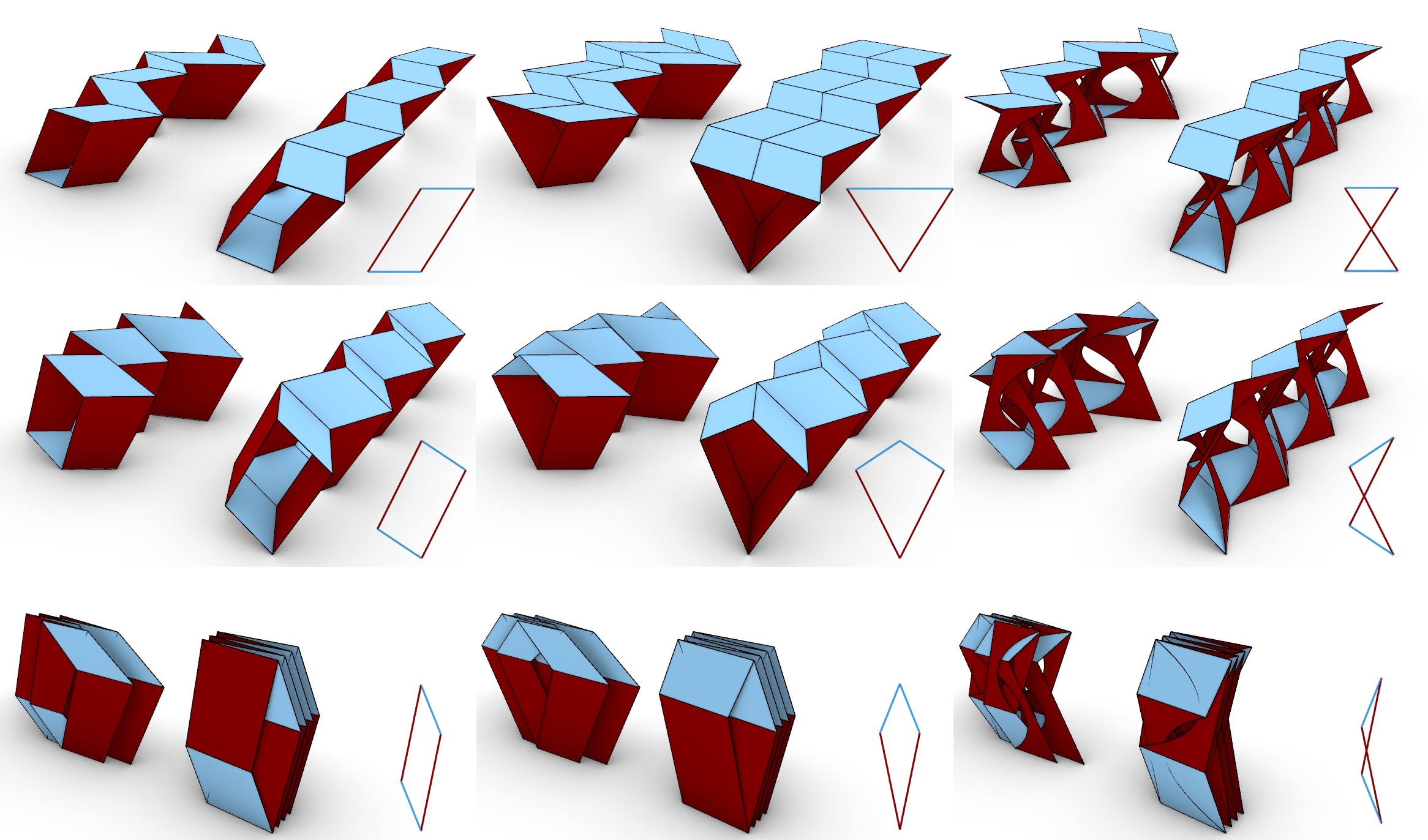}
{\small

}
\end{overpic}
    \caption{Illustration of T-hedral tubes with all types of possible quadrilateral cross-sections; parallelogram (left), deltoid (center) and anti-parallelogram (right). For the latter case we used special shaped 
    panels in order to avoid the self-intersection. For each cross-section, tubes are generated with regular and irregular zigzag trajectories and shown in three stages of folding, (from top to bottom) completely unfolded (up to the flexion limit), half-folded and flat-folded.   
    }
    \label{all-quad-profiles}

\end{figure}
This theorem also gives a pure geometric proof for the rigid-foldability of some tubes 
in the literature, where the continuous flexibility was only verified by numeric computations 
(cf.\  Appendix A of \cite{FPT2015}). 
Moreover, the characterization given in Theorem \ref{closed-discrete-profile} is more relaxed than the one of \cite{FPT2015}, which can easily be seen by the example of a deltoidal cross-section (Fig.\ \ref{all-quad-cross-sections}). This is the 
simplest closed polyline, which does not go along with the conditions given in \cite{FPT2015}.  Also, the anti-parallelogram as a cross-section is valid if it is 
symmetric with respect to the $X$-axis, but it yields a tube with self-intersection\footnote{The anti-parallelogram is contained in the solution set of \cite{FPT2015} but was not mentioned explicitly.}. 
In Fig.\ \ref{all-quad-cross-sections}, all types of possible quadrilateral cross-sections are illustrated and in Fig.\ \ref{all-quad-profiles}, examples of corresponding T-hedral tubes. 
An application is illustrated in Fig.\ \ref{Bench}, where a transformable bench is designed utilizing the properties of the flexion limits.

\begin{remark}
   Due to the conditions of Theorem \ref{closed-discrete-profile}, one can show that the line-segments of the cross-sections in each class, can always be divided into pairs with a same length and same absolute value of the slope.  Therefore, the cross-sections are polygons with even number of sides. However, polygons with odd number of sides (see \cite[Fig.\ 2]{Chen}) could be generated and can be seen as degenerated even-sided polygons with at least one interior angle equal to $\pi$. \hfill $\diamond$
\end{remark}

\begin{figure}[t]
\centering
\begin{overpic}[width=165mm]{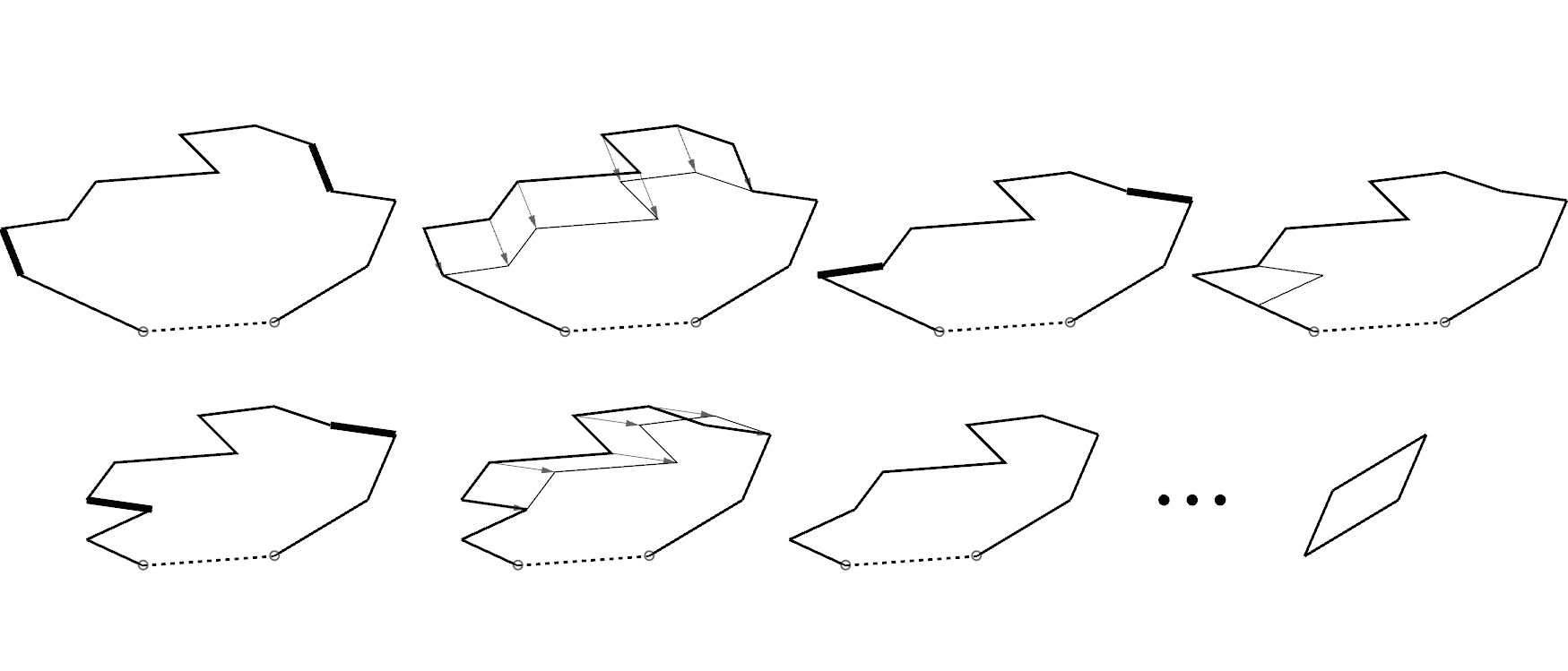}
{
\put(24,22){$\rightarrow$}
\put(50,22){$\rightarrow$}
\put(74,22){$\rightarrow$}
\put(74,22){$\rightarrow$}
\put(95,22){$\rightarrow$}
\put(2,8){$\rightarrow$}
\put(24,8){$\rightarrow$}
\put(46,8){$\rightarrow$}
\put(69,8){$\rightarrow$}
\put(80,8){$\rightarrow$}
\put(12,19){(a)}
\put(40,19){(b)}
\put(64,19){(c)}
\put(88,19){(d)}
\put(12,4){(e)}
\put(40,4){(f)}
\put(64,4){(g)}
\put(88,4){(h)}
}
\end{overpic}
 \caption{The process of cross-sectional reduction by removing pairs of segments. (a) An identical pair of segments is bold. (b) A row of parallelograms is being removed from the cross-sectional area. (c) A pair of vertical symmetric segments is bold. (d) A deltoid is being removed to construct an identical pair. (e,f,g) Repeating the same process as illustrated in (a,b,c) until we end up in (h). The remaining quad is always either a parallelogram or a deltoid or an anti-parallelogram. Not that, the dashed lines represent the part of the cross-section that has been removed for simplification.}
    \label{cross-section_tiling}
\end{figure}
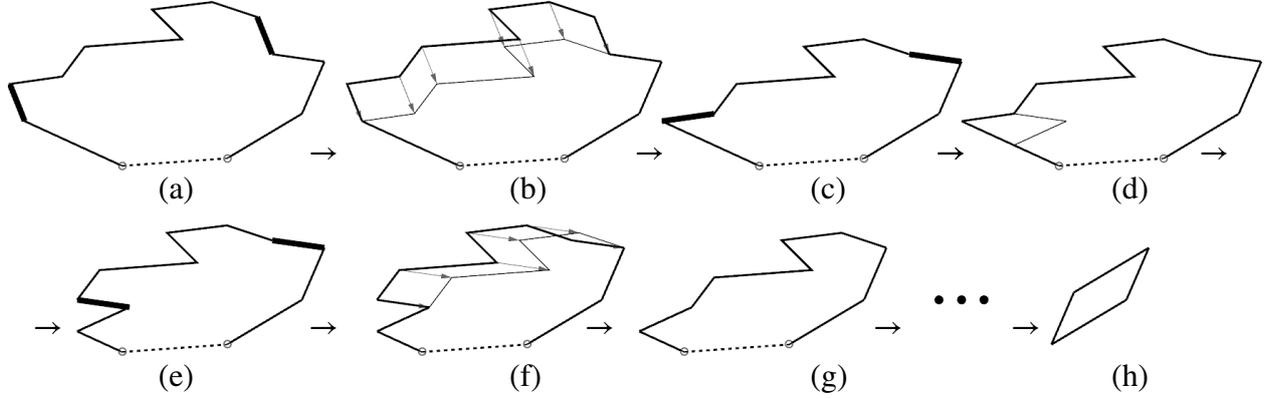


\begin{theorem}\label{profile-division}
A closed cross-sectional polyline satisfies the conditions of Theorem \ref{closed-discrete-profile} if and only if, it can be generated by a process of combining (adding and subtracting) parallelograms and deltoids, defined in \cite{Chen}.

\end{theorem}
\begin{proof}
Here we give a sketch of the proof by construction. For a given closed polyline, satisfying the conditions of Theorem \ref{closed-discrete-profile}, we can assume, the line-segments (up to division into smaller pieces) always come in pairs of identical segments or reflection symmetric with respect to the $Z$-axis. For each pair of identical segments (Fig. \ref{cross-section_tiling}(a)), we can cover the part of the boundary, enclosed between these two segments \footnote{ There are two boundary parts for each pair of segments, which are equally valid for this selection.}, with a row of parallelograms (Fig. \ref{cross-section_tiling}(b)). The result, is the new cross-section for the next iteration (cf.\ Fig.\ \ref{cross-section_tiling}(c)). The new cross-section, contains all the edges of the previous cross-section, only the pair of identical segments has been removed by a simple translation. For each pair of reflected segments (cf.\ Fig.\ \ref{cross-section_tiling}(c)), first we add (or subtract) a deltoid, to reflect it (cf.\ Fig.\ \ref{cross-section_tiling}(d)) and obtain a pair of identical segments (cf.\ Fig.\ \ref{cross-section_tiling}(e)). Then we repeat the process of removing as explained (cf.\ Fig.\ \ref{cross-section_tiling}(f,g)) until only a quad remains (cf.\ Fig.\ \ref{cross-section_tiling}(h)). It can be shown, that this quad is either a parallelogram, a deltoid, or an anti-parallelogram. Not that, the deltoids we use in this process must be oriented to have their symmetry line vertical to the base plane.

The reverse direction of the theorem, can be proven in a similar manner, starting from one of the three flexible quads of Fig.\ \ref{all-quad-cross-sections}, and showing that the process of adding quads of the same class, keeps the conditions given in Theorem \ref{closed-discrete-profile}. 
\end{proof}
\begin{figure}
\centering
\begin{overpic}[width=165mm]{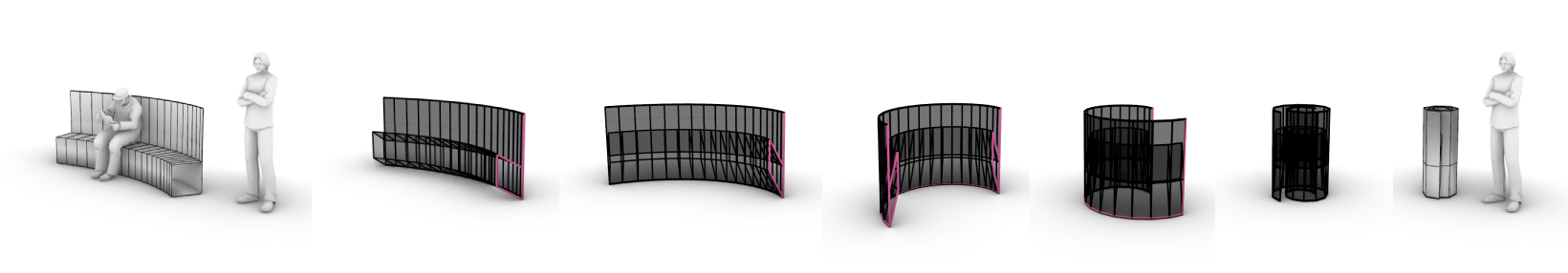}
\end{overpic}

 \caption{Design example: A parametric rigid-foldable spiral bench, which can be curled into a tight compact state.}
    \label{Bench}
\end{figure}

Theorem \ref{profile-division} shows that our method can generate all the cross-sections of \cite{Chen}.
This also opens doors for developing algorithms to approximate arbitrary given closed curves with flexible closed polylines. 

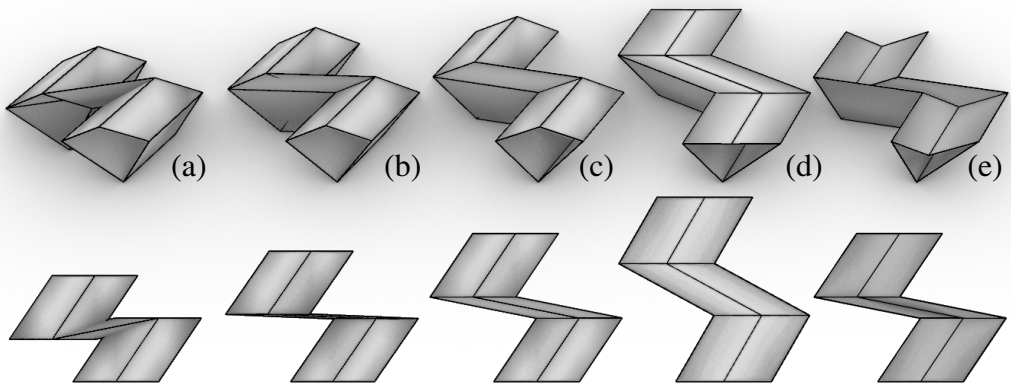
\begin{figure}[b]
\centering
\begin{overpic}[width=165mm]{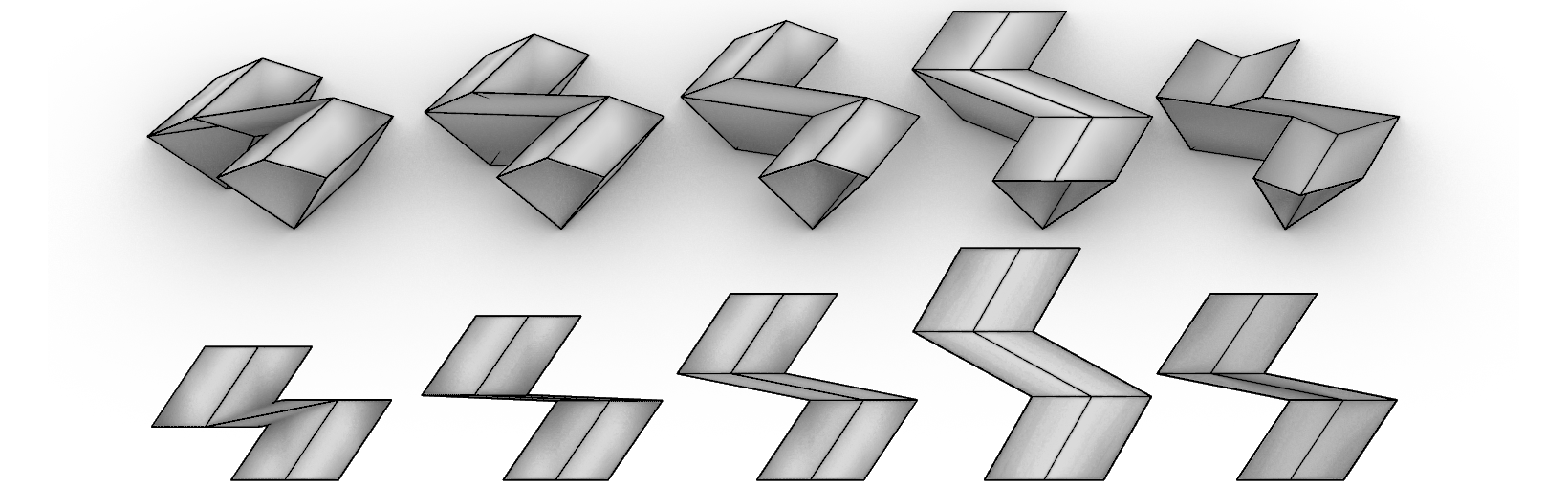}
{

 \put(22.5,17){(a)}
\put(39.5,17){(b)}
\put(55,17){(c)}
\put(71.5,17){(d)}
\put(86,17){(e)}

}
\end{overpic}
 \caption{A T-hedral tube in (c) the initial state, (b) its vertical and (d) horizontal flexion limits and in (a) vertical and (e) horizontal switched states. }
    \label{switches}
\end{figure}

\subsection{Switching}\label{subsec:switch}

As a consequence of Remark \ref{discrete_bifurcation} also T-heral tubes are in a bifurcation configuration, if either (1) a row is parallel to the base plane or  (2) a closed column, which we call loop, is completely flat (cf.\ Fig.\ \ref{switches}). The change of the working mode is also known as (1) horizontal and (2) vertical switching. 

\begin{itemize}
\item[ad 1)]
This kind of switching has to be done 
in accordance with the conditions stated in Theorem \ref{closed-discrete-profile}. Moreover, it was already used in \cite{FPT2015}, where it is also very well illustrated. 
\begin{remark}
If one is close to the horizontal bifurcation the structure can also snap from one mode into the other by minor deformations of the material. \hfill $\diamond$
\end{remark}
\item[ad 2)]
The vertical switching is only of theoretical nature as it causes self-intersections of the tube.
\end{itemize}

\section{Extension to the semi-discrete and smooth tubes with an isometric deformation}\label{sec:semitube}


\subsection{Smooth and semi-discrete analogs of T-hedra}\label{subsec:Smooth}

We replace the polylines, which define the geometry of a T-hedron (i.e.\ profile  polyline, trajectory polyline, prism polygon; cf.\ Section \ref{sec:T-hedra}), by corresponding smooth analogs (i.e.\ profile curve, trajectory curve, prism curve). As a consequence the discrete kinematic generation of the surfaces (cf.\ Section \ref{subsec:subclass}) gets continuous and we end up with a smooth surface, which can be classified into the same types as the T-hedra. In the literature the smooth analog of T-hedra are known as profile-affine surfaces \cite{sauer,graf} but we refer to them as {\it smooth T-surfaces} in the remainder of the paper for reasons of consistency in notation. According to \cite{sauer,graf,RT2021} each smooth T-surface allows a 1-parametric isometric deformation\footnote{
Beside this term we can still use the notion {\it continuous flexion} in the context of  smooth and semi-discrete T-surfaces, but we cannot name it  {\it rigid-foldability} anymore.} into smooth T-surfaces of the same type.  
Moreover, the flexion limits are reached \cite{sauer,graf} if and only if the smooth T-surface touches along a trajectory/profile curve the carrier plane of this curve (i.e.\ this carrier plane is the tangent plane along the trajectory/profile curve). 

If only the profile is replaced by a smooth curve, we end up with a so-called {\it semi-discrete T-surfaces of the vertical kind}, 
which consists of smooth cylindrical columns. 
A semi-discrete T-surface of the {\it horizontal kind} is obtained if only the trajectory and prism curve are smooth. Each resulting smooth row belongs to a ruled surface $\Lambda$ due to the straight line-segment of the profile polygon. In the general case these rulings envelope a spatial curve located on the guiding cylinder $\Gamma$, thus $\Lambda$ is a tangent surface. 
In special cases the ruled surface $\Lambda$ can degenerate to cones or cylinders. As a consequence each column-strip/row-strip of a semi-discrete T-surface of the vertical/horizontal kind is a developable strip (D-strip) and therefore these structures are of special interest for practical applications (cf.\ \cite{pottmann}). Moreover, the developability of the strips is an advantage for rapid prototyping of surface models (cf.\ \cite{maleczek2022rapid}). 

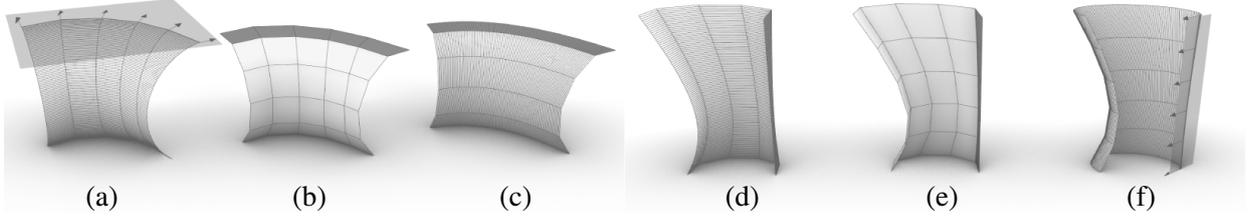
\begin{figure}
\centering
\begin{overpic}[width=165mm]{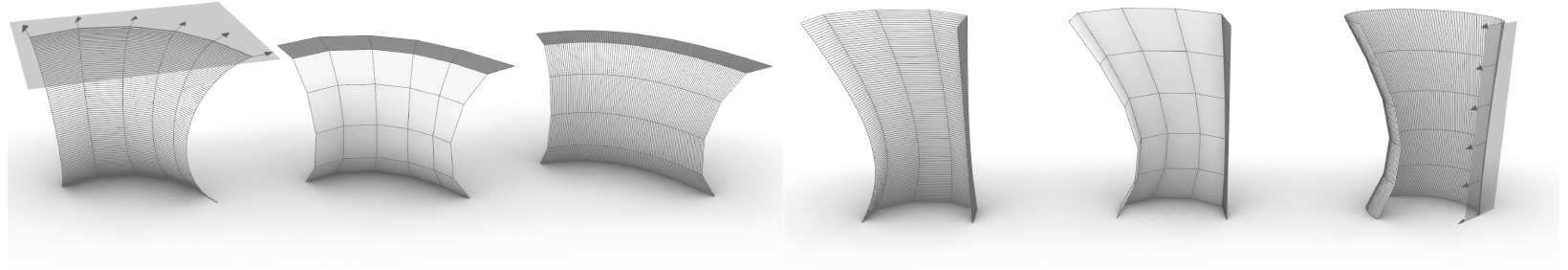}
{\small

\put(7,0.5){(a)} 
\put(23.5,0.5){(b)}
\put(40,0.5){(c)} 
\put(58,0.5){(d)} 
\put(74,0.5){(e)} 
\put(90,0.5){(f)} 
}
\end{overpic}
    \caption{Discrete and semi-discrete T-hedra  in horizontal (a,b,c) and vertical (d,e,f) flexion limits.}
    \label{flexion}

\end{figure}

The flexion limits for the semi-discrete T-surfaces can be argued from those of the discrete case on the one hand (see Section \ref{subsec:Iso}) and from those of the smooth case on the other hand noted above. 
Thus, for the semi-discrete version of the vertical kind, the
flexion limits are reached if either a D-strip gets completely flat (Fig.\ref{flexion} (d)) or a D-strip has a ruling with a horizontal tangent plane (Fig.\ref{flexion} (a)). For the horizontal kind, the flexion limits are reached if either a D-strip gets parallel to the base plane (Fig.\ref{flexion} (c)) or a D-strip has a ruling whose tangent plane touches along the entire associated profile polyline of the semi-discrete surface (Fig.\ref{flexion} (f)).

\begin{remark}
Note that Remark \ref{discrete_bifurcation} also holds true if one replaces PQ-strips by D-strips and T-hedron by semi-discrete T-surface, respectively.
\hfill $\diamond$
\end{remark}

\subsection{Curved cross-section}\label{subsec:Curve}
Using the cross-sections of Theorem \ref{closed-discrete-profile} as profile polygons together with smooth trajectory and prism curves result in the construction of semi-discrete T-tubes of the horizontal kind (cf.\ Fig.\ref{semi}(c)). 

For smooth T-tubes and semi-discrete T-tubes of the vertical kind one has to smoothen the cross-section. In this case  the 
profile curve cannot be smooth everywhere during the isometric deformation 
as the existence of a horizontal tangent corresponds to a flexion limit. 
Therefore  we can assume without loss of generality that there exists
a state of deformation, where the profile has no horizontal tangent.  
\begin{lemma}\label{lem:twoangle}
We suppose that a closed cross-section $c$ consists of two 
$C^1$ continuous curve segments $c_1$ and $c_2$. Moreover the existence of vertical 
tangents is only allowed in the points $K_1$ and $K_2$, where the two curves are 
stitched together. Then $c$ is compatible with the isometric deformations of T-type if and only if one of the following four cases hold true (up to translations): 
\begin{enumerate}[i.]
\item
$c_1$ and $c_2$ are identical ($c_2$ = $c_1^i$),
\item
$c_1$ and $c_2$ are reflection symmetric with respect to the $X$-axis ($c_2$ = $c_1^{ii}$),
\item
$c_1$ and $c_2$ are reflection symmetric with respect to the $Z$-axis ($c_2$ = $c_1^{iii}$),
\item
$c_1$ and $c_2$ are related by a $180^{\circ}$ rotation ($c_2$ = $c_1^{iv}$).
\end{enumerate}

\end{lemma}

\begin{figure}
\centering
\begin{overpic}[width=140mm]{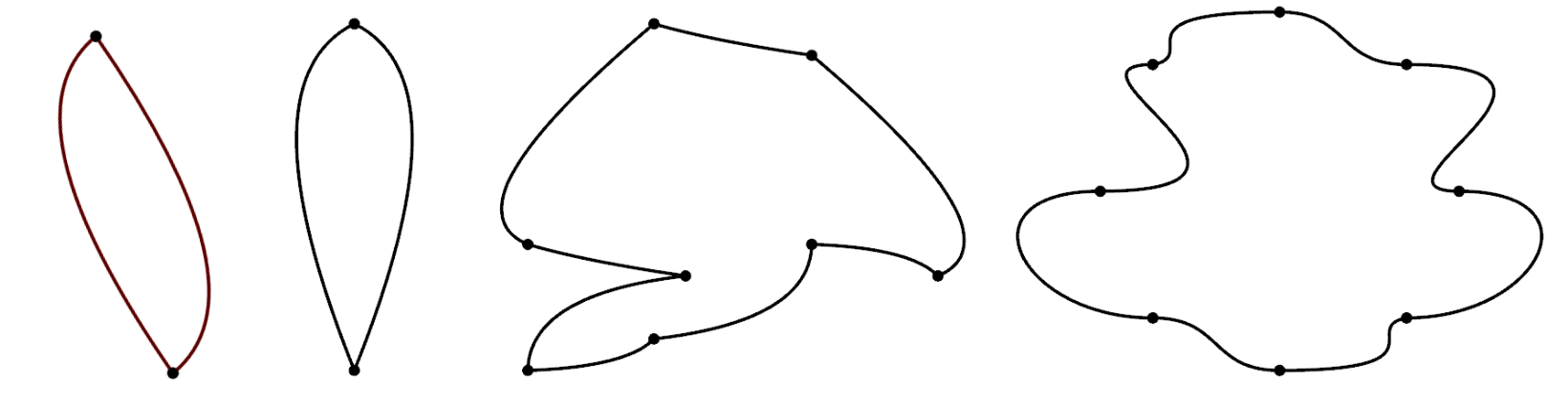}
{
\put(3,10){$c_1$}
\put(11,16){$c_1^{iv}$}
\put(17,10){$c_1$}
\put(26,10){$c_1^{iii}$}
\put(34,21){$c_1$}
\put(39,10){$c_2$}
\put(32,6){$c_3$}
\put(38,1){$c_4$}
\put(46,24.5){$c_2^i$}
\put(58,17){$c_1^{iii}$}
\put(55,6.5){$c_4^{ii}$}
\put(49,4){$c_3^{iv}$}
\put(74,26){$c_1$}
\put(72,16){$c_2$}
\put(68,7.5){$c_3$}
\put(75,2){$c_4$}
\put(85,25){$c_4^{i}$}
\put(94.5,16){$c_2^{ii}$}
\put(96,5){$c_3^{iii}$}
\put(87,0){$c_1^{iv}$}
\put(7,0){$K_1$}
\put(19,0){$K_1$}
\put(6,24.5){$K_2$}
\put(22,25){$K_2$}
}
{\small
\put(8,-4){(a)} 
\put(21,-4){(b)}
\put(44,-4){(c)} 
\put(80,-4){(d)} 

}
\end{overpic}
\newline
    \caption{
    Cross-sections (profiles) consisting of $C^1$ continuous curve-segments, which are compatible with the isometric deformation of T-type. (a,b) Simplest cross-sections with minimum number of curve-segments. (c,d) Cross-sections with all four cases of Lemma \ref{lem:twoangle}.}
    \label{All-curves&Smooth-profile}

\end{figure}
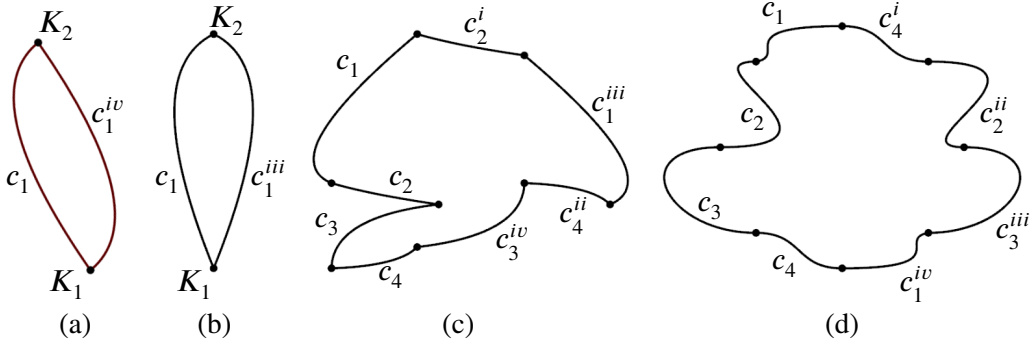

\begin{proof}
Without loss of generality we can assume an arc-length parametrization $(f_i(s),0,g_i(s))$  with arc-length $s\in[0;k_i]$ of the curve-segment $c_i$ for $i=1,2$ with 
$c_i(0)=K_1$ and $c_i(k_i)=K_2$  (cf.\ Fig.\ \ref{All-curves&Smooth-profile}(a,b)). It is further no restriction to assume that $K_1$ is located within the origin of the coordinate frame.
According to \cite{RT2021} the curve $c_i$ deforms as follows with respect to the 
transformation parameter $t$:
\begin{equation}\label{eq:deform}
\left(tf_i,0,\int_0^{s}\sqrt{(1-t^2)f_i'^2+g_i'^2}ds\right)
\end{equation}
where the prime is indicating the derivative with respect to $s$. 
Now for all $t$ the $Z$-coordinates at $s=k_i$ have to coincide; i.e.
\begin{equation}\label{eq:equalz}
\int_0^{k_1}\sqrt{(1-t^2)f_1'^2+g_1'^2}ds = \int_0^{k_2}\sqrt{(1-t^2)f_2'^2+g_2'^2}ds
\end{equation}
From the special value $t=0$ it can be seen that the curve-segments $c_1$ and $c_2$ unroll along the $Z$-axis, which shows that 
the arc lengths of $c_1$ and $c_2$ have to match; i.e. $k:=k_1=k_2$. Thus from Eq.\ (\ref{eq:equalz}) 
the following condition results after differentiation with respect to $s$ and squaring:
\begin{equation}
(1-t^2)(f_1'^2-f_2'^2)+(g_1'^2-g_2'^2)=0. 
\end{equation}
This equation can only be fulfilled independent of $t$ for $f_1'=\pm f_2'$ and $g_1'=\pm g_2'$. 
This already shows that the curve-segments $c_1$ and $c_2$ have to meet one of the four relations listed in the lemma, as there cannot be a switch in the signs due to the assumptions of $C^1$ continuity and non-existence of horizontal tangents and 
vertical ones (which are only possible in the endpoints of the curve segments).
\end{proof}

We continue with a definition of considered cross-sections and their admissible decompositions. 
\begin{definition}\label{def:decomposition}
We consider a cross-sections $c$, which is a $C^0$ continuous finite compositions of at least $C^1$ smooth curve segments $c_1,\ldots, c_n$. 
An admissible decomposition of $c$ is obtained by introducing a finite set of break points, 
which contains at least all points with either (a) a $C^1$ discontinuity 
or (b) a vertical tangent. If the curve $c$ contains a vertical line-segment $c_i$ 
then both end-points are assigned as break points but not the points in between. 
\end{definition}


\begin{remark}
Note that points having a vertical tangent keep this property under the isometric deformation.  Moreover we can assume that in all states of the deformation with exception of the flexion limits (cf.\ Fig.\ \ref{All-curves&Smooth-profile}(d)),there do not exist points with a horizontal tangent (cf.\ Fig.\ \ref{All-curves&Smooth-profile}(c)). In these states the local maxima and minima of the cross-section with respect to the $Z$-coordinate have to be singular points ($C^1$ discontinuity).
\hfill $\diamond$ 
\end{remark}
Based on Lemma \ref{lem:twoangle} and Definition \ref{def:decomposition} we can give the following theorem: 

\begin{theorem}\label{thm:general}
A closed cross-section $c$ of Definition  \ref{def:decomposition} is compatible with the isometric deformations of T-type if there exists an admissible decomposition such that the following conditions hold: 

\begin{enumerate}
\item
Curve-segments which are related in one of the four possibilities stated in (i-iv) belong to the same class.
\item
For each class we add up the signum values of the individual curve-segments $c_i$, where $\sign(c_i)$ is defined as follows: We start in an arbitrary break point and run counter-clockwise through the cross-section. 
Curve-segments $c_i$, which are strictly monotonically increasing (resp.\ decreasing) during the sweep, result in $\sign(c_i):=+1$ (resp.\  $\sign(c_i):=-1$). 
\item
For each class this sum has to be zero.
\end{enumerate}
\end{theorem}
\begin{remark}
Note that this general formulation also includes the discrete setting stated in Theorem \ref{closed-discrete-profile}.
\hfill $\diamond$ 
\end{remark}

By combining the cross-sections of Theorem \ref{thm:general} with smooth trajectory and prism curves, we obtain smooth T-tubes (cf.\ Fig.\ref{semi}(a)), which have at least two horizontal creases (singular curves). 
If the trajectory and the prism are polygonal then we obtain semi-discrete T-tubes of the vertical kind (cf.\ Fig.\ref{semi}(b)).

\begin{remark}
Horizontal switching is only possible if one of the curve segments  $c_1,\ldots, c_n$ 
mentioned in Definition \ref{def:decomposition} is straight (and becomes horizontal during the flexion). 
The vertical switching is still faced with 
the problem mentioned in
Section \ref{subsec:switch}.
\hfill $\diamond$

\end{remark}
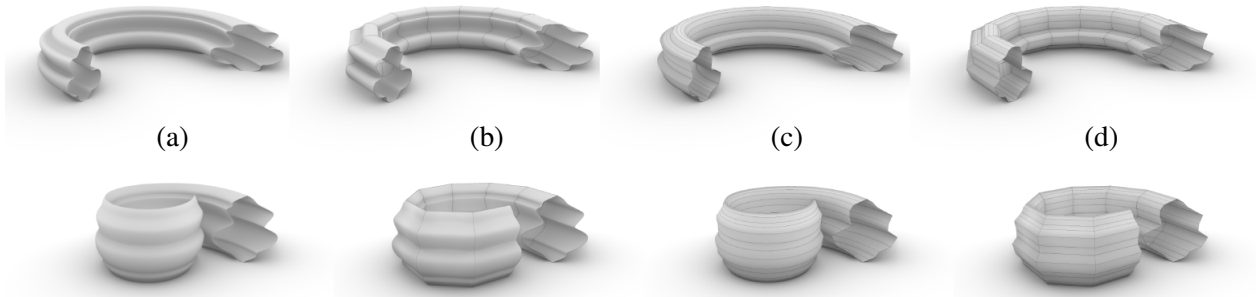
\begin{figure}
\centering
\begin{overpic}[width=165mm]{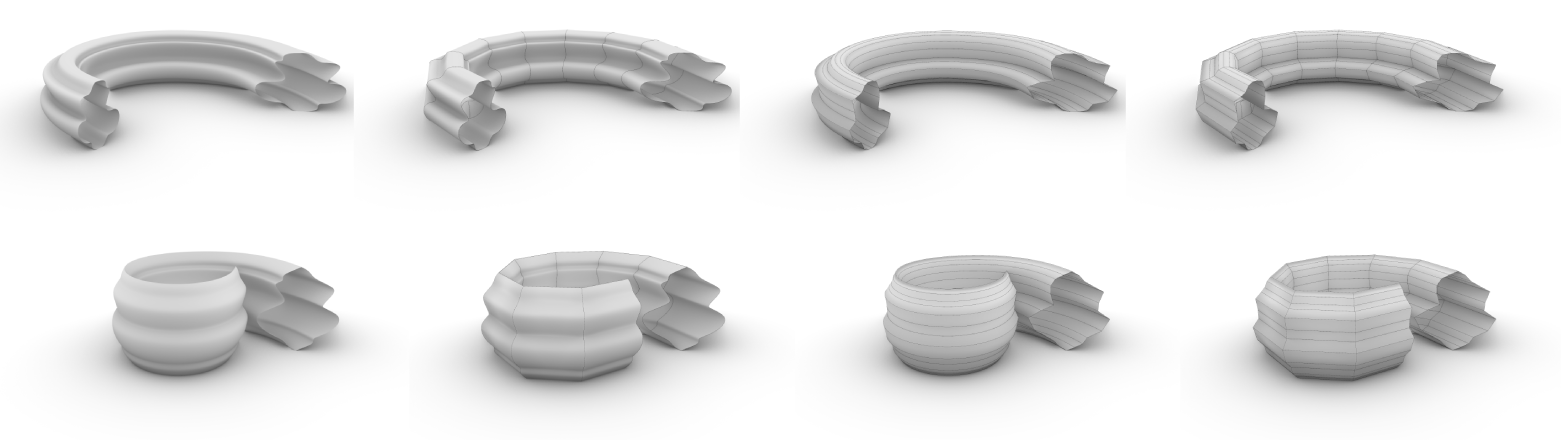}
{\small
\put(12,15){(a)} 
\put(37,15){(b)}
\put(61,15){(c)} 
\put(86,15){(d)} 

}
\end{overpic}
 \caption{ A smooth tube and its discretizations in initial (top row) and deformed (bottom row) states. (a) Smooth. (b) Semi-discrete of vertical kind. (c) Semi-discrete of horizontal kind. (d) Discrete.}
    \label{semi}
\end{figure}


\section{Structures composed of tubes}\label{sec:Structure}

In the literature not only rigid-foldable tubes are reported but there are also some 
publications dealing with their composition to surfaces and metameterials in a way that the
property of being continuous flexible is preserved. 
In the following we review the existing tubular surfaces and metamaterials with respect to the property of being T-hedral.

\begin{figure}
\centering
\begin{overpic}[scale=0.78]{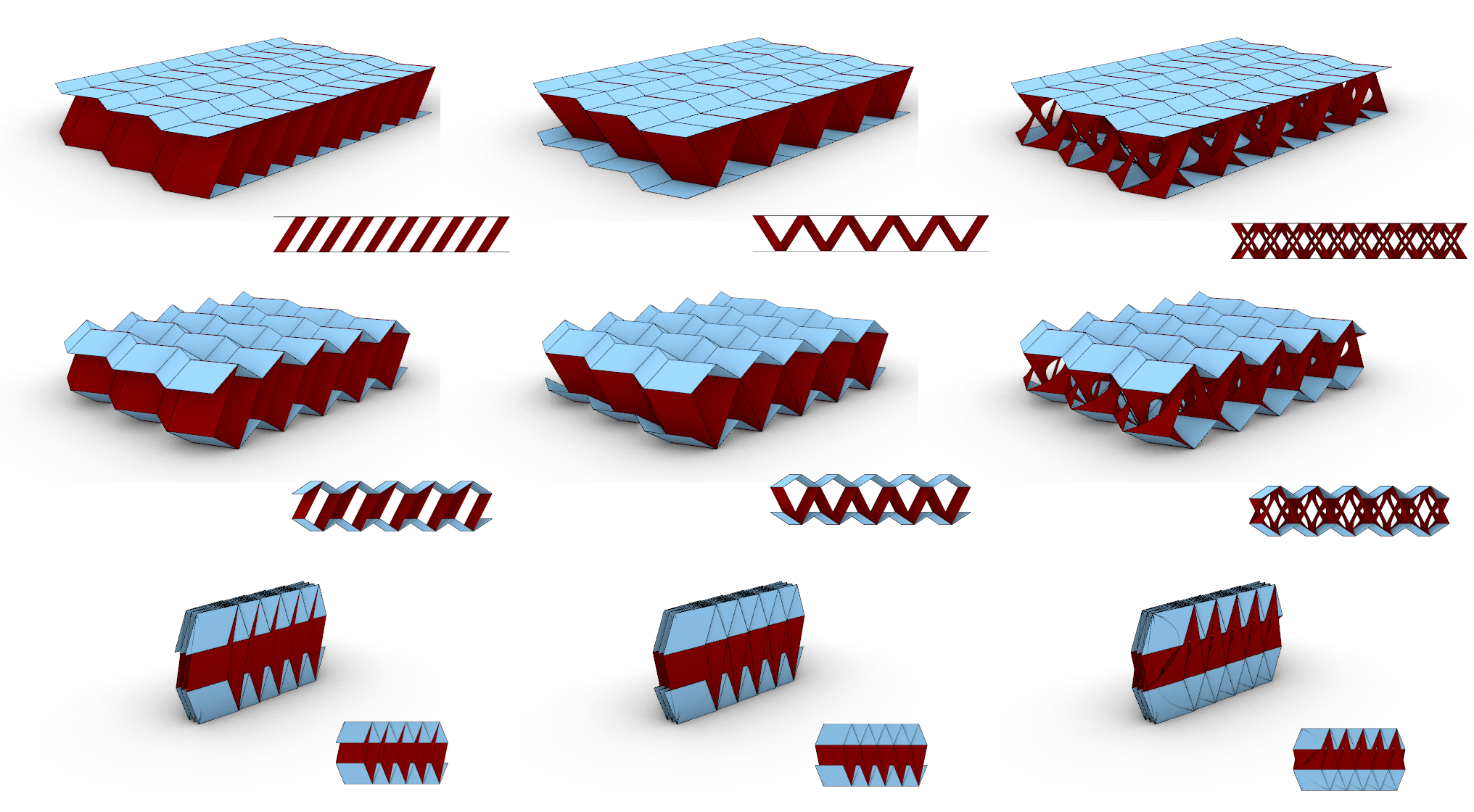}
{\small

}
\end{overpic}
    \caption{Sandwich structures made of aligned-coupled flat-foldable T-hedral tubes with all possible quad cross-sections in tree stages of folding, (from top to bottom) completely unfolded (up to the flexion limit), half-folded and flat-folded. The sandwich structure in the middle (with deltoid cross-sections) also appears in \cite{guest2013}.}
    \label{all-sandwiches}

\end{figure}

\subsection{Surfaces composed of tubes}\label{sec:surf}
First we discuss surfaces, where one can distinguish two different compositions, namely those which share edges or faces, respectively:

\begin{enumerate}
\item
Edge-sharing surfaces: The authors of \cite{ishizawa} distinguish the {\it Miura-ori surface} and the {\it eggbox surface}. If the surface approximate a plane then both structures are T-hedral.
If the surface approximate a vault (cf.\ Section 4 of \cite{ishizawa}) then only the eggbox structure is T-hedral. The Miura-ori vault is not T-hedral and can be folded rigidly 
due to the reasons given in Footnote \ref{octahedron}. 
\item
Face-sharing surfaces: One can distinguish the following two types of coupling of tubes forming the surface. 
    \begin{enumerate}
    \item Aligned-coupling:
    Most examples of this type given in the literature are T-hedral; e.g. the planar case as well as the vault case of \cite{ishizawa}, the surfaces given by Tachi in \cite[Fig.\ 11]{Tachi2009} and 
    \cite[Fig.\ 11]{tachi_iass2012}, the two-layered structures given in \cite{klett_iass2016} and  
    the tubular arches of \cite{gattasetal2017}.
    See also the examples given in Fig.\ \ref{all-sandwiches}. 
    
    However, there are also surfaces of aligned-coupled tubes, which are not T-hedral, like the top or bottom layer of the sandwich structure illustrated in \cite[Fig.\ 7]{tachietal2015}. Until now this is the only example known to the authors, but it is again composed of the special tubes mentioned in Footnote \ref{octahedron} which also reasons the rigid-foldability of the structure.
    \item
    Zipper-coupling:
    The complete zipper structure is not T-hedral, but in most of the 
    examples given in the literature \cite{pnas} it consists of T-hedral tubes. 
    The only  exception known to the authors, are again the zipper coupling of tubes with a
    parallelogram cross-section mentioned in Footnote \ref{octahedron} (cf.\ \cite{Filipovetal2019} and 
    \cite{tachietal2015}).
    This observation motivates us to study the zipper coupling of T-hedral tubes in more detail in Section \ref{sec:Zipper}.
    \end{enumerate}
\end{enumerate}

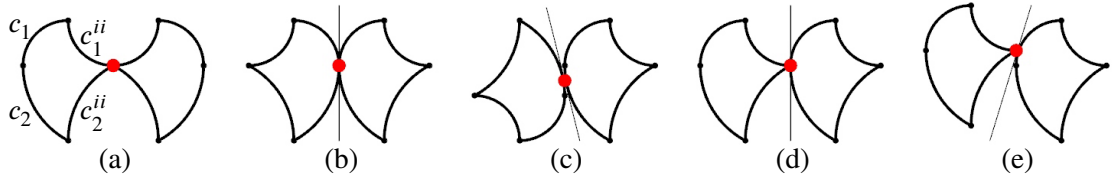
\begin{figure}
\centering
\begin{overpic}[width=165mm]{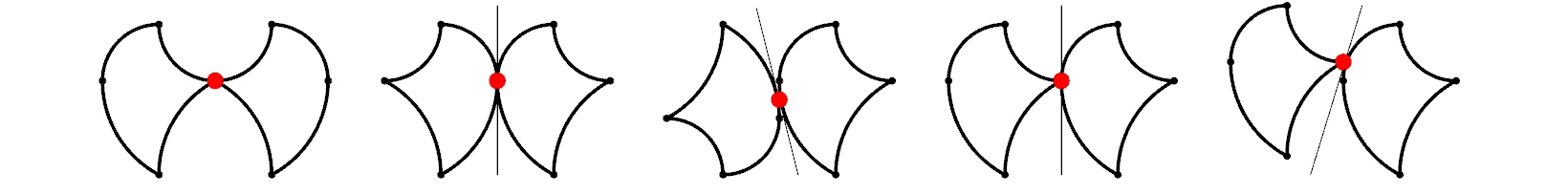}
{\small
\put(5.3,10){$c_{1}$} 
\put(5.3,3){$c_{2}$} 
\put(11,9){$c_{1}^{ii}$}  
\put(11,3){$c_{2}^{ii}$}   
\put(12.5,-1){(a)}
\put(30.5,-1){(b)}
\put(48.5,-1){(c)}
\put(66.5,-1){(d)}
\put(84.5,-1){(e)}
}
\end{overpic}
    \caption{All possible types of generalized edge-sharing. (a) Crease-sharing. (b) Curves touching on points where they share a vertical tangent which is not necessarily their symmetry line, and (c) sharing a non-vertical tangent. (d) A crease touching a curve on a point with a vertical tangent, and (e) with a non-vertical tangent. 
    }
    \label{Edge-couplings}

\end{figure}
\begin{figure}
\centering
\begin{overpic}[width=165mm]{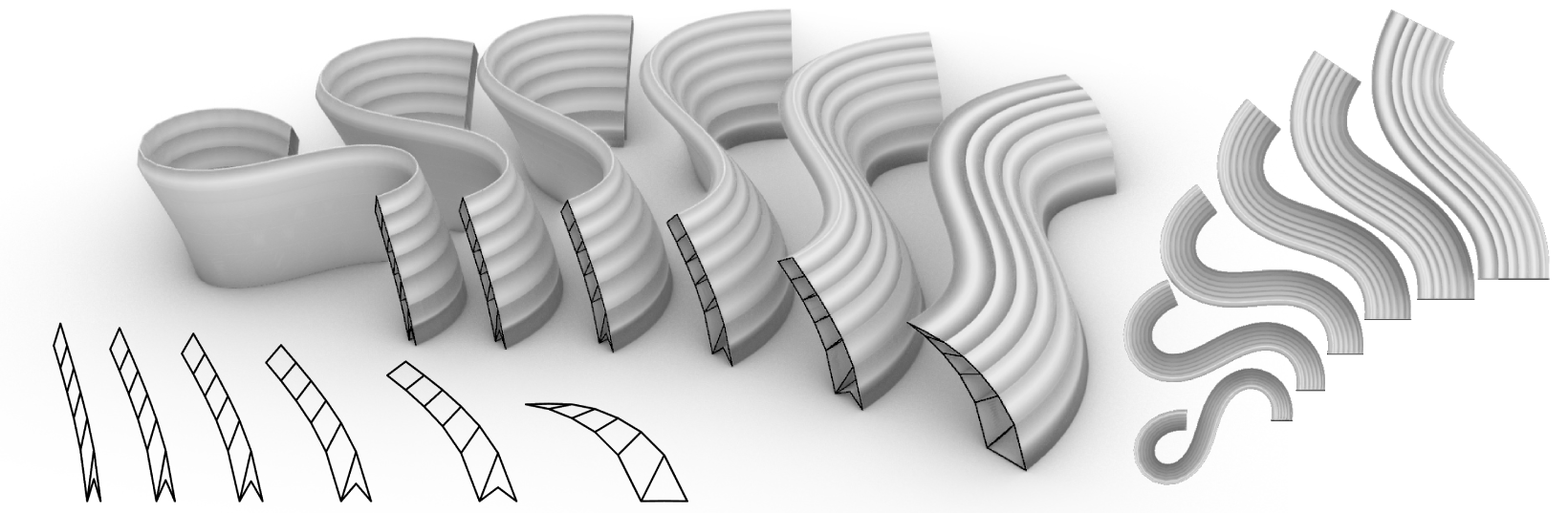}
{\small

}
\end{overpic}
    \caption{Folding sequence of a 
    surface composed of semi-discrete T-tubes of the horizontal kind sharing faces by an aligned-coupling in front, perspective and top views (from left to right).}
    \label{Sandwich-pavilion}

\end{figure}

The identification of tubular surfaces to be T-hedral either of the edge-sharing type or face-sharing type with aligned-coupling already implies a pure geometric reasoning of their rigid-foldability, which was in some cases only argued on base of numerical computations 
(cf.\ \cite{ishizawa}). 
Moreover, this geometric insight enlarge the design space of edge-sharing and face-sharing surfaces with an aligned-coupling considerably, and in addition it can be generalized to the semi-discrete and 
smooth case in the following way:
\begin{enumerate}
\item
{\bf Generalized edge-sharing surface:} The extension of edge-sharing surface construction to semi-discrete tubes 
of the horizontal kind is straight forward. For the vertical kind and the smooth case one can distinguish different types, which result from the different kinds of trajectory curves; namely ($\alpha$) singular curves (creases) and regular curves with ($\beta$) vertical tangent-planes with ($\gamma$) non-vertical tangent-planes. 
If two tubes share a common trajectory curve, then the tubes also have to touch along that curve in 
order to avoid local penetration of the tubes. This implies the following possible combinations:
    \begin{enumerate}[$\star$]
        \item ($\alpha$)-($\alpha$): We can call this a crease-sharing  (see Fig.\ \ref{Edge-couplings}(a)).
        \item ($\beta$)-($\beta$): The tangent-planes along that touching curves of two tubes are vertical  (see Fig.\ \ref{Edge-couplings}(b)), as this property is kept under the deformation.
        \item ($\gamma$)-($\gamma$): The touching curves are the trajectories of corresponding points of curve segments $c_1$ and $c_2$ which are related by a $180^{\circ}$ rotation (relation iv of Section \ref{sec:semitube}). An example for this is illustrated  in Fig.\ \ref{Edge-couplings}(c).
        \item ($\alpha$)-($\beta$): In this case the common curve is a crease of one tube and a regular curve with vertical tangent-planes of the other tube (see Fig.\ \ref{Edge-couplings}(d)).
        \item ($\alpha$)-($\gamma$): In this case the common curve is a crease of one tube and a regular curve with non-vertical tangent-planes of the other tube (see Fig.\ \ref{Edge-couplings}(e)).
    \end{enumerate}

\item
{\bf Generalized face-sharing surfaces with aligned-coupling:} The extension to the smooth case and to both kinds of the semi-discrete case are straight forward and illustrated in Fig.\ref{Sandwich-pavilion}.
Note that a 
semi-discrete surface of the horizontal kind was 
reported in \cite[Fig. 7]{klett_iass2016}.
\end{enumerate}

\begin{remark}
Note that all the structures listed in Sections \ref{sec:surf}, 
which are rigid-foldable due to their composition 
of tubes with parallelogram cross-section mentioned in Footnote \ref{octahedron}, have the potential to 
be generalized according to the aleady mentioned idea illustrated in \cite[Fig.\ 06]{tachi_transformables_2013}. \hfill $\diamond$
\end{remark}

\subsection{Metamaterials composed of tubes}\label{sec:meta}

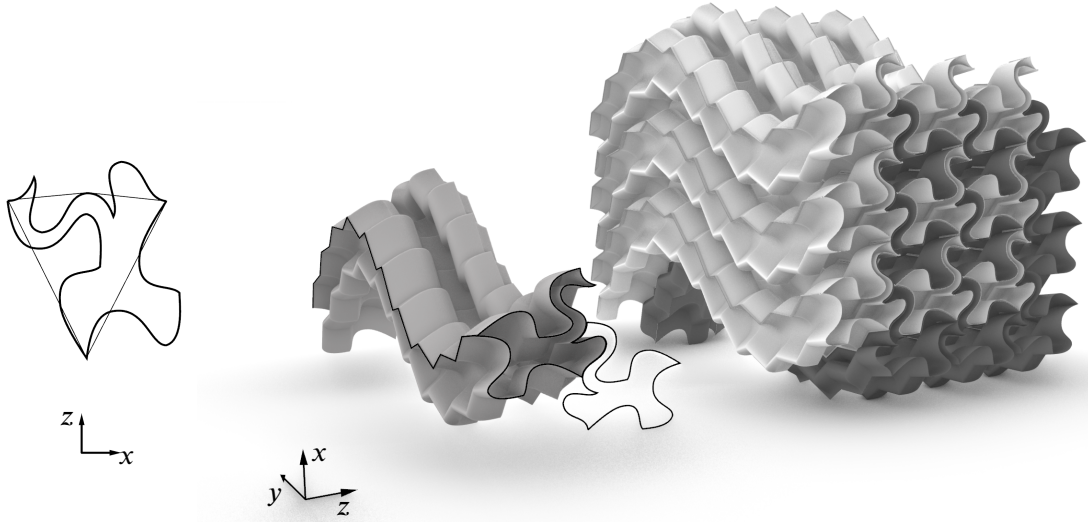
\begin{figure}
\begin{center}
\begin{overpic}
[width=25mm]{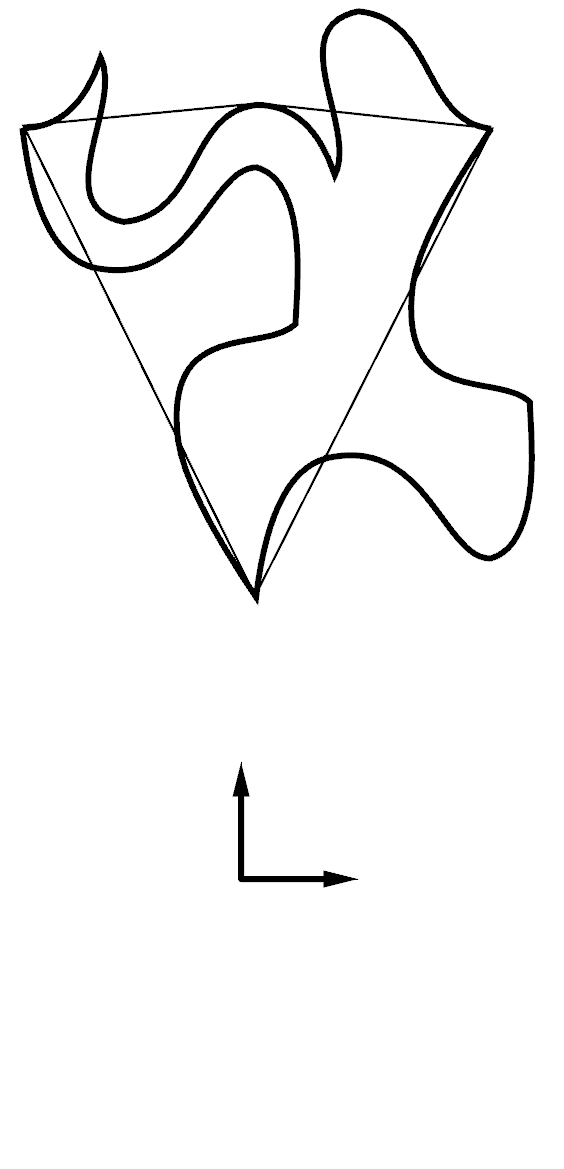}
\begin{small}
\put(30,22){$x$}
\put(15,32){$z$}
 \end{small}
\end{overpic}
\begin{overpic}[width=130mm]{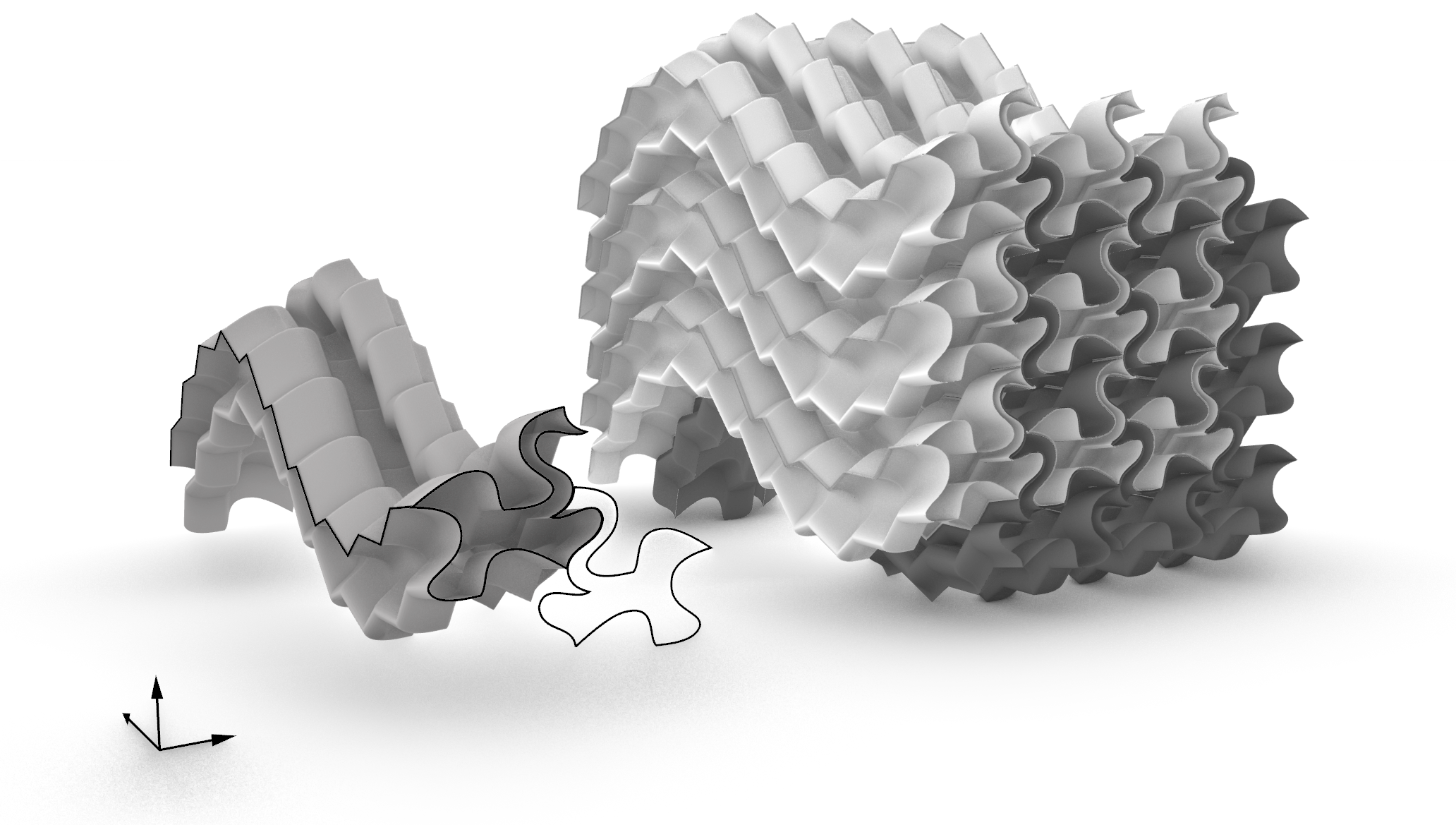}
\begin{small}
\put(14,3.5){$z$}
\put(11.5,9){$x$}
\put(7,5){$y$}
\end{small}
\end{overpic}
\end{center}
 \caption{Semi-discrete T-hedral tubes of the vertical kind, which are
 aligned-coupled to a flat-foldable metamaterial.
  (left) a cross-section on a deltoid base inspired by M.C.\ Escher's Moebius birds (center) a single tube and (right) the resulting metamaterial.}
    \label{Escher}
\end{figure}

We use the following classification for the review of matamaterials composed 
of tubes:
\begin{enumerate}
    \item 
    Aligned-coupling: Under this term we summarize edge- and face-sharing structures, as this separation can depend on the interpretation; e.g.\ the structure given in Fig.\ \ref{Escher} can be seen as an edge-sharing one by only looking at the white tubes, but also as a face-sharing structure by considering black and white tubes. 
    Examples for aligned-coupling are  \cite[Figs.\ 10 and 24]{tachi_iass2012}, the multi-layer materials
    of \cite[Figs.\ 10 and 11]{klett_iass2016}, the Miura-ori matamaterial of \cite{guest2013}, the modular foldable structures of \cite{davood}.
    All these examples are T-hedra, thus the T-hedral construction unifies all 
    these approaches and enlarges the design space of these metamaterials. 
    Clearly, the generalization to the semi-discrete cases and smooth case is again straight forward (cf.\ Fig.\ \ref{Escher}).
    \item
    Woven-coupling: The structure of \cite[Fig.\  24]{tachi_iass2012}, which was obtained by an aligned-coupling, can also be decomposed in a different way; namely into two families of T-hedral tubes. Neighboring tubes of the same family are edge-sharing and  neighboring tubes of different families are face-sharing (cf.\ \cite[Fig.\  24]{tachi_iass2012}). No further examples of 
    tubular woven metamaterials are known to the authors so far.
    \item
  Zipper-coupling: The tubes of the zipper metamaterial in  \cite[Fig. 5]{pnas} are T-hedral as well as in \cite[Fig. 7]{Filipovetal2019}. For example sandwich surfaces in \cite[Fig. 7]{tachietal2015} and \cite[Fig. 8(a,c)]{Filipovetal2019} could be extended to non-T-hedral zipper metamaterials by adding further layers.
    The generalization of the zipper coupling of T-hedral tubes is discussed in detail in the next section. 
\end{enumerate}

\section{T-hedral zipper tubes}\label{sec:Zipper}

The term {\it zipper tubes} always refer to a pair of discrete tubes (in most examples both tubes are congruent) consisting of planar quads where each interior vertex has valence four. 
Each of these tubes has a one-parametric isometric deformation. Moreover the tubes can be glued together along a row, which is called 
{\it zip row} in the remainder of the paper, and the resulting structure is still rigid-foldable. 

As already pointed out in the last section, most of the examples given in the literature consist of T-hedral tubes. 
This perspective allows us to  explain the underlying ``zipper principle'' in more detail and to give a  precise definition for 
so-called {\it T-hedral zipper tubes}, which reads as follows:

\begin{definition}\label{def:zip} 
{\it T-hedral zipper tubes} consist of two T-hedral tubes, which can be glued together along the zip row such that the resulting structure is still isometric deformable. Moreover, the base planes of the two involved T-tubes are not parallel/identical. 
\end{definition}

\begin{remark}
If the base planes of the two involved T-tubes are parallel/identical then one ends up with the 
already mentioned face-sharing aligned-coupling (cf.\ Section \ref{sec:surf}). \hfill $\diamond$ 
\end{remark}

Note that Definition \ref{def:zip} is written in a way 
that it also holds for zipper structures composed of 
smooth or semi-discrete T-tubes. In the following we show 
that smooth zipper structures cannot exist, by proving the
following theorem.

\begin{theorem}\label{thm:profil}
The profile curve of a zip row has to be a straight line-segment.
\end{theorem}

\begin{proof}

If the profile is curved then it spans a unique carrier plane which equals the profile plane. The intersection of (infinitesimal) neighboring 
profile planes gives the axis of the (infinitesimal) stretch rotation, which has to orthogonal to both non-parallel/identical base planes 
mentioned in Definition \ref{def:zip}. 
This is only possible if this axis is at infinity, thus we have an 
(instantaneous) translation causing a T-hedron of type I.
In the following we will show that even this is impossible. 

Without loss of generality we can assume again (cf.\ Section \ref{subsec:Curve}) that the curved profile of a zip row has
an arc-length parametrization  $\Vkt p(s):=(f(s),0,g(s))$.  Then this profile deforms (deformation parameter $t$) isometrically 
into 
\begin{equation}\label{eq:p_deform}
\Vkt p_t(s):=\left(tf,0,\int_0^{s}\sqrt{(1-t^2)f'^2+g'^2}ds\right)
\end{equation} 
according to \cite{RT2021}, where  the prime is indicating the derivative with respect to $s$.
Note that the $XY$-plane is the base plane of this deformation. 

Instead of computing a deformation of $\Vkt p(s)$ with respect to another base plane, we rotate 
the profile curve about the $Y$-axis by an arbitrary angle $\alpha\in]0;\pi[ \cup ]\pi;2\pi[$. 
In this way we get $\mVkt p(s):=(\cos{\alpha}f(s)+\sin{\alpha}g(s),0,\cos{\alpha}g(s)-\sin{\alpha}f(s))$, 
which can again be deformed (deformation parameter $\dach t$) in analogy to Eq.\ (\ref{eq:p_deform}) which yields:
\begin{equation}\label{eq:p_deform2}
\mVkt p_t(s):=
\left(\dach t(\cos{\alpha}f+\sin{\alpha}g),0,\int_0^{s}\sqrt{(1-\dach t^2)(\cos{\alpha}f'+\sin{\alpha}g')^2
+(\cos{\alpha}g'-\sin{\alpha}f')^2}ds\right)
\end{equation} 
Due to the fundamental theorem of planar curves we only have to check if the curvature functions $\kappa_t(s)$ and $\dach\kappa_t(s)$ 
of the curves $\Vkt p_t(s)$ and $\mVkt p_t(s)$ are identical. As these curves have still arc-length parametrizations we 
can use the well-known formula:
\begin{equation}
\kappa_t(s)^2=\| \Vkt p_t(s)''\|^2, \quad
\dach\kappa_t(s)^2=\| \mVkt p_t(s)''\|^2, 
\end{equation}
which implies the condition $\| \Vkt p_t(s)''\|^2-\| \mVkt p_t(s)''\|^2=0$. 
Straight forward computation of this condition yields an expression of the form $K_1^2K_2=0$ with
\begin{align}
K_1:&=f'g'' - f''g', \\
K_2:&=\dach t^2(t^2 - 1) (\sin{\alpha} f' -  \cos{\alpha}g')^2 + t^2(1-\dach t^2)g'^2.
\end{align}
Therefore either $K_1=0$ has to hold, which is discussed later on, or $K_2$ vanishes. In the latter case we can express $\dach t^2$ which yields
\begin{equation}
\dach t^2=\frac{t^2g'^2}{g'^2-(t^2 - 1) (\sin{\alpha} f' -  \cos{\alpha}g')^2}.
\end{equation}
As $\dach t$ has to be constant for all values $s$, the partial derivative of $\dach t^2$ with respect to $s$ has to vanish (for all values of $t$). 
This implies the condition:
\begin{equation}
t^2(t^2-1)\sin{\alpha}g'(g'\cos{\alpha}-f'\sin{\alpha})K_1=0.
\end{equation}
As $\alpha\in]0;\pi[ \cup ]\pi;2\pi[$ has to hold we only remain with the following three possibilities:
\begin{enumerate}
\item
$g'=0$ implies that $g$ is constant. Geometrically the profile curve is $\Vkt p(s)$ is a 
straight line parallel to the $X$-axis. 
\item
$g'\cos{\alpha}-f'\sin{\alpha}=0$ is fulfilled for with $f'=\cos{\alpha}$ and $g'=\sin{\alpha}$ as the side-condition $f'^2+g'^2=1$ of 
arc-length parametrization has to hold. Integration shows that the profile curve is a straight line with inclination $\alpha$. 
\item
The system of differential equations 
$K_1=0$ and $f'^2+g'^2=1$ can also be solved which yields the solution 
\begin{equation}
f=ks+d,\quad g=+\sqrt{1-k^2}s+e
\end{equation}
with $k,d,e\in\RR$, which also expresses a straight line.  
\end{enumerate}
This finishes the proof of Theorem \ref{thm:profil}.
\end{proof}

In view of Definition \ref{def:zip} we have to compute zip rows with the property that the two planar rims remain planar during an isometric deformation of the zip row. 
Due to Theorem \ref{thm:profil} one can restrict to the case that the zip row belongs either to a cone (type II and type III; cf.\ Subsection \ref{sec:relation}) or to a cylinder (type I), respectively. 
This geometric problem was discussed in detail by the last author in \cite{N2022} and its results are summarized in the following theorem:

\begin{theorem}\label{thm:results}
In the following cases there exists an isometric deformation $\iota$ of a cylinder/cone such that the planar cuts of two non-parallel planes $\alpha$ and $\beta$ remain both planar:
\begin{enumerate}
    \item 
    In the cylindrical case there is no condition to the shape of the cylinder, but the intersection line of the two planes $\alpha$ and $\beta$ has to be orthogonal to the ruling direction of the cylinder. This holds for the smooth and discrete case. 
    \item
    In the discrete conical case the cone has to be a cap of a Bricard octahedron of the plane-symmetric type. Moreover, the planes $\alpha$ and $\beta$ have to be orthogonal to the plane of symmetry $\omega$ of this cone (cf.\ Fig.\ \ref{non-translational zipper} (left)). 
\end{enumerate}
\end{theorem}
\begin{figure}
\centering
\begin{overpic}[width=165mm]{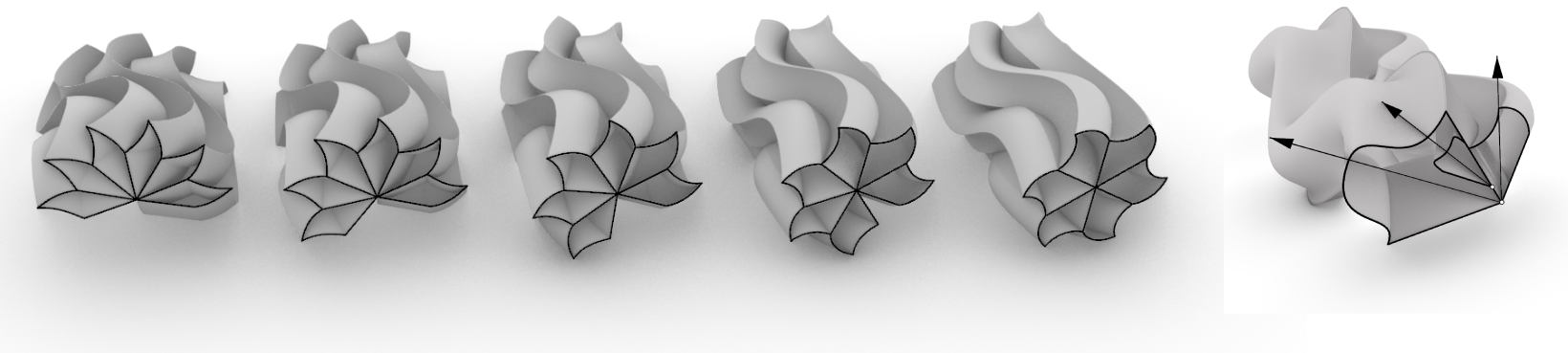}
{
\put(38,2){$(a)$}
\put(89,2){$(b)$}
}
\end{overpic}

 \caption{ 
 (a) The sequence of isometric deformation of zipper-coupled tubes with semi-discrete zip rows. The outer rows of the tubes can also be seen as smooth surfaces, as they have been created by a curved segmented profile and smooth trajectory. (b) Three tubes sharing a single zip row, having three different base planes with illustrated normal vectors. }
    \label{semi-multiple-zipper}
\end{figure}
\begin{remark}\label{rem:zipper}
For the smooth conical case we solved the resulting system of partial differential equations symbolically up to a final ordinary differential equations, which we were not able to solve \cite{N2022}. \hfill $\diamond$
\end{remark}

 Moreover it can easily be seen that all planes belonging to the pencil spanned by $\alpha$ and $\beta$ (and all planes parallel to them) 
 cut the cone/cylinder in curves which also remain planar during the deformation. Therefore one can, not only zip two tubes along a zip row but an arbitrary number (Fig.\ref{semi-multiple-zipper} (b)).

\begin{remark}\label{rk:glue}
Note that along a zip row one can also glue two T-hedral surfaces with different base planes. For the translational case this is not exciting as then one just fall out of type I into type II or III. But for the non-translatinoal case 
this is a non-trivial coupling and it is not clear yet to which class of \cite{Ivan} the corresponding $(3\times 3)$ complex with the zip row in the middle belongs.
\hfill $\diamond$
\end{remark}

Based on Theorem \ref{thm:results} we are able to present generalized T-hedral zipper tubes in the following two subsections.


\subsection{Generalization of T-hedral zipper tubes of the translational type }\label{subsec:Translational-zipper}

Until now only T-hedral zipper tubes of the translational type appeared in 
the literature \cite{Filipovetal2019,pnas}, where the trajectory curve is a regular zigzag and the profile curve is a rhombus. But the zipper coupling can easily be generalized as follows: 

One generates a 
translational zip row according to item 1 of Theorem \ref{thm:results} and constructs two closed profiles which meet the conditions of Theorem  \ref{closed-discrete-profile} or Theorem \ref{thm:general} and contain an overlapping part of the line-segment (cf.\ Theorem \ref{thm:profil}) of the zip row. Note that this construction also implies the possibility of semi-discrete zip rows (see Fig. \ref{semi-multiple-zipper} (a)).
Moreover, Theorem \ref{thm:results} also implies a proof for the flexion of the resulting zipper tubes, which was only reasoned by numerical computations in \cite{Filipovetal2019,pnas} to the best knowledge of the authors.


\begin{figure}[t]
\begin{center}
\begin{overpic}
    [height=55mm]{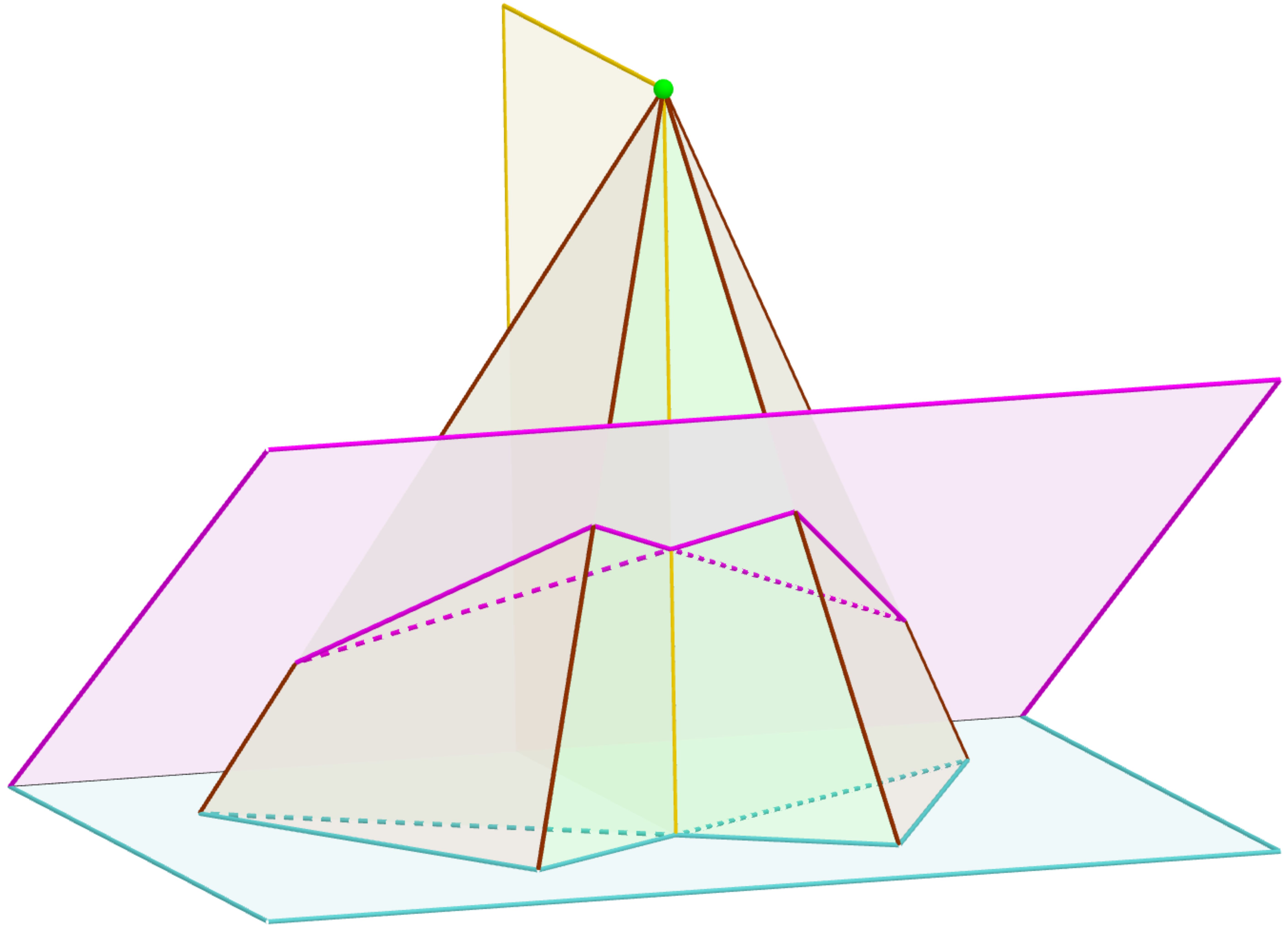}
\begin{small}
\put(53,64){$V$}
\put(40,65.5){$\omega$}
\put(92,38.5){$\alpha$}
\put(85,7){$\beta$}
\put(24.5,3){$b_i$}
\put(32.5,29){$a_i$}
\put(16,22){$A_i$}
\put(10,7){$B_i$}
\put(31.5,20){$a_{i+1}$}
\put(29,11){$b_{i+1}$}
\end{small}     
  \end{overpic} 
\quad
\begin{overpic}
    [height=60mm]{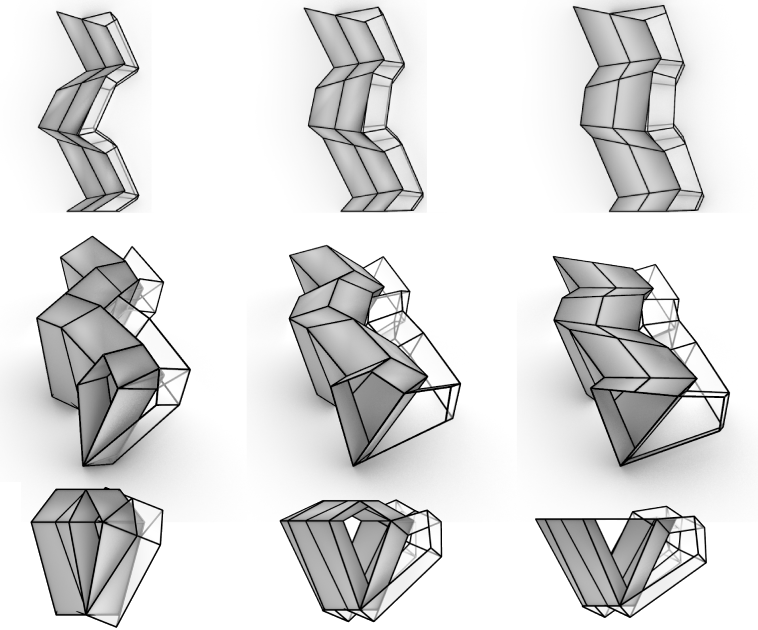}
  \end{overpic} 
\end{center}	
\caption{
(left) Cap of a Bricard octahedron of the plane-symmetric type. (right) 
Folding sequence of T-hedral zipper tubes of the non-translational type in top, perspective and front views (from top to bottom).}
  \label{non-translational zipper}
\end{figure}
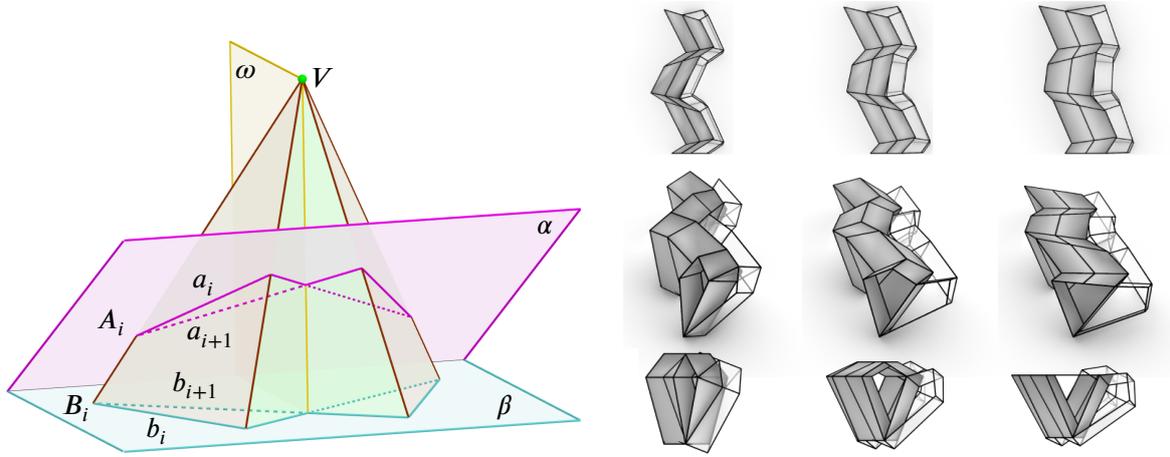    


\subsection{T-hedral zipper tubes of the non-translational type}\label{subsec:Non-translational-zipper}

For the first time we present T-hedral zipper tubes of the non-translational type, which are zipped along a discrete zip row. Their construction is done as follows:

According to \cite{N2022} the faces $f_1,\ldots,f_n$ of a discrete zip row of type II wrap multiple times (depending on $n$) the cap of a Bricard octahedron of the plane-symmetric type with vertex $V$, such that the rims $a_i$ and $b_i$ are located in the planes $\alpha$ and $\beta$, respectively, which are orthogonal to the plane of symmetry $\omega$. The line-segments $a_i$ and $a_{i+1}$ (resp.\ $b_i$ and $b_{i+1}$) are linked by the vertex $A_i$ (resp.\ $B_i$) for $i=1,\ldots, n-1$ (cf.\ Fig.\ \ref{non-translational zipper} (left)).

\begin{remark}\label{rmk:firstlast}
Note that the first face $f_1$  (resp.\ last face $f_n$) can be cut arbitrarily but we can assume that this is done in a way that the 
resulting start (resp.\ end) points $A_0$ and $B_0$ (resp.\  $A_n$ and $B_n$) are collinear with $V$. \hfill $\diamond$
\end{remark}

Therefore the points $A_1,\ldots, A_{n-1}$ with 
$A_j=A_{j+4}=A_{j+8}=\ldots$ for $j=1,\ldots,4$ are the vertices 
of an anti-parallelogram which is passed multiple times (depending on $n$).
The same holds true for the vertices $B_1,\ldots, B_{n-1}$.

The profile plane $\pi_{a_i}$ (resp.\ $\pi_{b_i}$)  through $a_i$ (resp.\ $b_i$) for the tube with base plane $\alpha$ (resp.\ $\beta$) is obtained by the plane orthogonal to $\alpha$ (resp.\ $\beta$) containing the line-segment $a_ib_i$. 

\begin{remark}
From this construction it is clear that the planes $\pi_{a_i}$ and $\pi_{b_i}$ only coincide if $a_i$ and $b_i$ are located in $\omega$. 
The construction of the zipper tubes becomes simpler if the two profiles of the tubes are within the same plane, thus we assume $A_0$ and $B_0$ to 
be located in $\omega$ (which is no restriction according to Remark \ref{rmk:firstlast}). \hfill $\diamond$
\end{remark}

The design space of the resulting zip row seems to be very narrow, but one must not forget the possibility of applying the map $\phi^{-1}$ (cf.\ Section \ref{sec:relation}), which turns the zip row from type II to type III. This enlarges the space of possible zip row designs considerably (cf.\ Fig.\ \ref{non-translational zipper} (right)).

\begin{remark}
Note that a cap of a  Bricard octahedron of the plane-symmetric type has two flat poses, which imply two flat poses of the corresponding zip row (of type II or type III). They are flexion limits as well as bifurcation configurations of the zipper tubes according to Section \ref{subsec:Iso}. \hfill $\diamond$
\end{remark}

\begin{figure}[b]
\centering
\begin{overpic}[width=165mm]{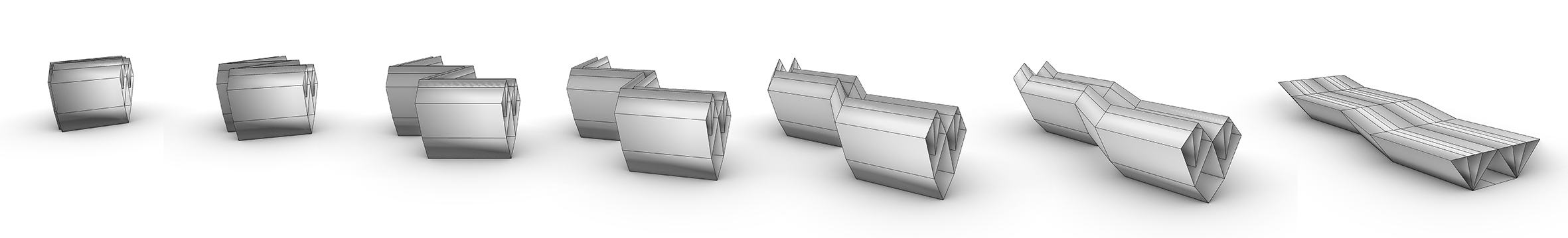}

\end{overpic}
    \caption{ Folding sequence of a T-hedral aligned-coupled tubular bridge with a nested deltoidal cross-section from the flat-folded state to the completely unfolded one (from left to right).
    }
    \label{bridge-folding}

\end{figure}

\section{Potential applications and outlook
}\label{sec:bridge}

\begin{figure}
\centering
\begin{overpic}[scale=0.5]{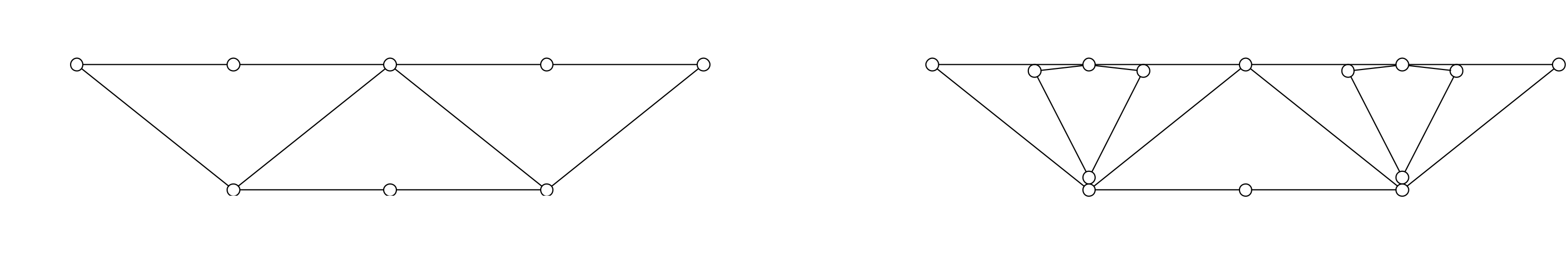}

\end{overpic}
    \caption{ Cross-section with standard base shape (left) and with inner deltoids to prevent bifurcation (right). 
    }
    \label{bridge sandwich}

\end{figure}

Among other potential use cases, the authors investigated different versions of  flat-foldable bridges consisting of  T-hedral tubes with a flat drive- or pathway in its final unfolded  state (cf.\ Fig.\ \ref{bridge-folding}).  The advantage of bringing an entire bridge or parts of it in flat-folded or developed configuration on site and unfold it there to its final shape, reduces a lot of grey energy and is therefore a widespread strategy, whereas transportation size is a key parameter. 
For this particular investigation the authors focus was the structural performance based on the global geometry of the bridge that should span 20 meters. The actuation or exact detailing is the topic of ongoing investigation. As mentioned above, T-hedral  structures have due to their geometric constraints, bifurcational behaviour, when faces are parallel to the XY-Plane. As this is the plane, which cars or pedestrians would use, and bring load and therefore actuation force to the mechanism, the authors developed a nested deltoid cross-section, that would block the structure in final state, and prevent from bifurcation  (cf.\ Fig.\ \ref{bridge sandwich}).

\begin{figure}
\centering
\begin{overpic}[scale=0.5]{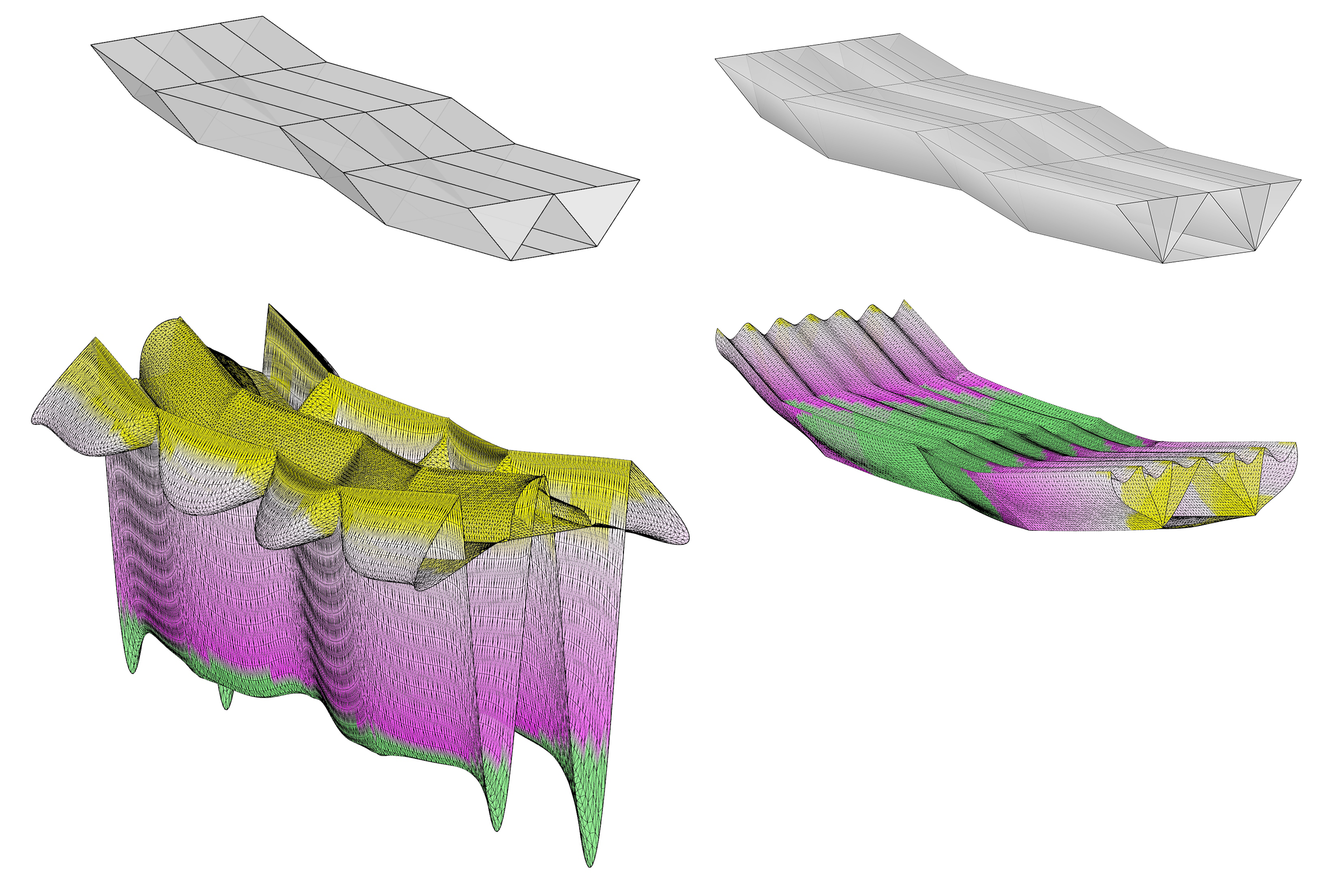}

\end{overpic}
    \caption{FEM Simulation of the two different cross-sections, showing the better performance through inner deltoid configuration (right) in comparison to standard shape (left)
    }
    \label{Bridge_Karamba}

\end{figure}

\begin{figure}
\centering
\begin{overpic}[width=165mm]{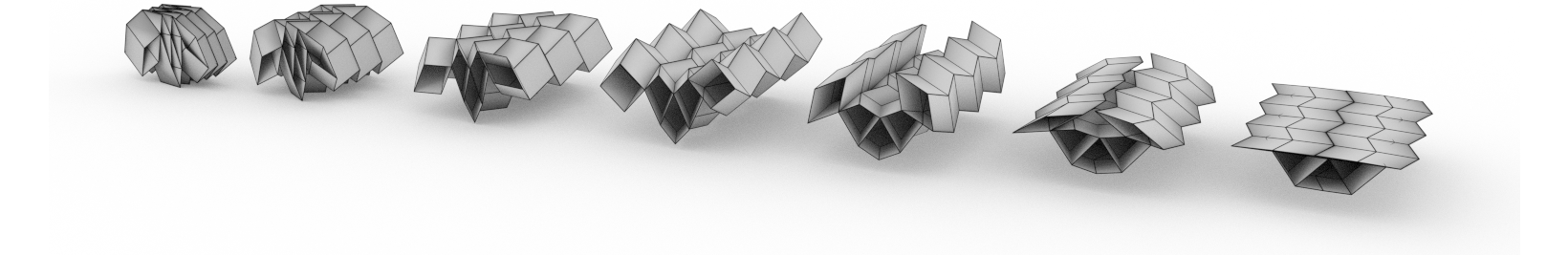}

\end{overpic}
    \caption{ Folding sequence of a T-hedral  aligned- and zipper-coupled tubular bridge from the flat-folded state to the completely unfolded one (from left to right).  
    The two zippered tubes on the side increase stiffness during motion according to \cite{pnas}, and prevent bifurcation as well as expand the flat pathway in the unfolded state.
    }
    \label{zipper-bridge}

\end{figure}

To have an insight in the global geometric stability of this nested configuration a parametric model in Rhino and Grasshopper was developed and investigated with the FEM-Tool Karamba \cite[Fig.\  24]{ Preisinger2013}. Assuming that both bridge beams in (cf.\ Fig.\ \ref{Bridge_Karamba}), have the same span, identical height and identical weight, the inner deltoid configuration has a better performance than the standard cross-section. 
The parametric setup developed by the authors will allow for a more detailed investigation to all bridge-beam parameters as  e.g. number of undulations, its shift in X-Axis, the inner deltoid-configuration or even  zipper tubes  variations  (cf.\ Fig.\ \ref{zipper-bridge}).   The ongoing research investigates if the advantages of flat folded transportation on site, might justify the reduction of stiffness of the final configuration.

Further future research directions associated with T-hedral tubes are given in the next subsection.
\subsection{Future research} \label{sec:future}
\begin{figure}[t]
\centering
\begin{overpic}[width=165mm]{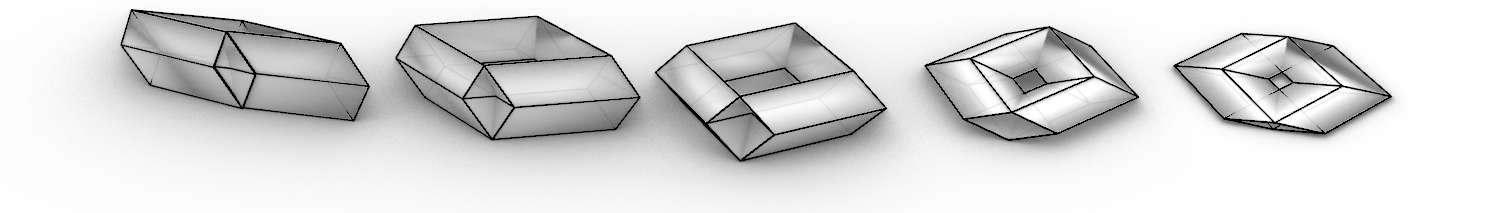}
\end{overpic}
    \caption{Folding sequence of the simplest closed T-hedral tube with a rhombic profile and trajectory. This flexible polyhedron is homeomorphic to a torus and also appears in \cite[Fig.\ 30.5]{pak}.  
    }
    \label{Torus}

\end{figure}

In the review part (cf.\ Section \ref{subsec:Review}) we have mentioned that the tubes of Miura and Tachi \cite{miura_synthesis} differ fundamentally from the presented construction method, but on the other hand we have not included them into the list of exceptions given in Section  
\ref{sec:general}. 
The reason for this is that these tubes are also T-hedral, but in contrast to all other T-hedral tubes mentioned so far in this paper not the profile curve is closed but the trajectory curve. A deeper study of this type of T-hedral tubes and the related surfaces and matamaterials (e.g.\ \cite[Figs.\ 20 and 31]{tachi_iass2012}) is dedicated to future research. We only want to mention that this allows us to
\begin{enumerate}
    \item 
    construct isometrically deformable closed tubes (see Fig.\ \ref{Torus}) which have self-intersections. To the best knowledge of the authors, the only example of this kind reported in the literature so far is the  ``Renault style'' polyhedron presented by Pak in \cite{pak}.
    \item
    come up with coupled tubes; by using e.g.\
    \begin{enumerate}
    \item 
    different profile curves along the same trajectory resulting in a kind of multi-layer tube, 
    \item
    one profile along different closed trajectories.
    \end{enumerate}
\end{enumerate}
Beside these issues we also want to study in more detail woven structures and surfaces, respectively, composed of T-hedral tubes. 
These surfaces, as well as the bridge-like structures, will be part of a more detailed investigation concerning the further development of building components and/or structures. 
In this context we also want to clarify the following questions: 
How does the angle enclosed by the base planes of the T-hedral zipper tubes (cf.\ Definition \ref{def:zip}) influence the stiffness of the system reported in \cite{pnas}.
From the geometric point of view the best performance is expected for a right angle. 
Moreover, can the stiffness be increased by the usage of non-translational zipper tubes as they are more interlaced than translational ones? 

\subsection{Open problems}
\begin{enumerate}
\item 
We conjecture that the word  ``{\it if}'' in the Theorems \ref{closed-discrete-profile} and \ref{thm:general}  can be replaced by ``{\it if and only if}''.
In the following we give some considerations to this problem with respect to
Theorem \ref{closed-discrete-profile}:

Each of the oriented line-segments $c_i$ for $i=1,\ldots,n$ can be seen as a vector
\begin{equation}
(s_i\cos{\alpha_i},0,s_i\sin{\alpha_i})    
\end{equation}
which deforms according to the formula given in Eq. (\ref{eq:deform}) by:
\begin{equation}
\left(ts_i\cos{\alpha_i},0,\sign(\sin{\alpha_i})s_i\sqrt{1-t^2{\sin{\alpha_i}}^2}\right)
\end{equation}
Therefore the following two conditions have to be fulfilled for all $t$:
\begin{equation}\label{eq:explain}
t\sum_{i=1}^ns_i \cos{\alpha_i}=0, \qquad 
\sum_{i=1}^n \sign(\sin{\alpha_i})s_i\sqrt{1-t^2{\sin{\alpha_i}}^2} =0.
\end{equation}
In contrast to the left condition the right one is not very nice due to the 
appearance of the square root expression. For $n=4$ all solutions (parallelogram, anti-parallelogram and deltoid) are known according to \cite{N2015}. But for the case  $n>4$ this remains totally open. 

The problem is similar for Theorem \ref{thm:general} but even more complex as we end up with a system of differential equations.
\item
Are there further classes of isometrically deformable non-cylindrical tube structures (with quadrilateral faces in the discrete case) beside those mentioned in the paper? Exist tube structures with non-planar quadrilateral faces?
Moreover, it is an open question if isometrically deformable closed tubes without 
self-intersection exist (cf.\ Subsection \ref{sec:future}).
\item
The open question mentioned in Remark \ref{rk:glue}.
\item 
Are there further zipper couplings beside those mentioned in the paper?
\item
Due to Remark \ref{rem:zipper} it remains open if non-translational semi-discrete zip rows of the horizontal kind exist.
\end{enumerate}

\noindent
{\bf Acknowledgement}
 The research is supported by grant F77 (SFB “Advanced Computational Design”) of the Austrian Science Fund FWF. 
\bibliography{tubes}

\end{document}